\documentclass{article}

\usepackage[numbers]{natbib}

\usepackage[final]{neurips_2021}

\usepackage{verbatim}
\usepackage[utf8]{inputenc} 
\usepackage[T1]{fontenc}    
\usepackage{hyperref}       
\usepackage{url}            
\usepackage{booktabs}       
\usepackage{amsfonts}       
\usepackage{nicefrac}       
\usepackage{microtype}      
\usepackage{xcolor}         
\usepackage{amsmath}
\usepackage{amssymb}
\usepackage[colorinlistoftodos,bordercolor=orange,backgroundcolor=orange!20,linecolor=orange,textsize=scriptsize]{todonotes}
\usepackage{enumitem}
\usepackage{pifont}
\usepackage{algorithm,algorithmicx,algpseudocode}
\usepackage{subcaption}
\usepackage{wrapfig}

\def\R{\mathbb{R}}

\def\R{\mathbb R}

\newcommand{\cO}{\mathcal{O}}

\newcommand{\norm}[1]{\left\| #1 \right\|}

\newcommand{\cbraces}[1]{\left( #1 \right)}

\newcommand{\braces}[1]{\left\{ #1 \right\}}
\newcommand{\EndProof}{\begin{flushright}$\square$\end{flushright}}

\def\<#1,#2>{\langle #1,#2\rangle}
\newcommand{\sqn}[1]{\norm{#1}^2}
\newcommand{\one}{\mathbf{1}}

\newtheorem{lemma}{Lemma}
\newtheorem{theorem}{Theorem}
\newtheorem{definition}{Definition}

\newtheorem{assumption}{Assumption}
\newtheorem{corollary}{Corollary}
 
 \title{Distributed Saddle-Point Problems Under  Similarity
}
\author{%
  Aleksandr Beznosikov \\
  MIPT\thanks{Moscow Institute of Physics and Technology}, HSE University and Yandex, Russia\\
  \texttt{anbeznosikov@gmail.com} \\
  \And
  Gesualdo Scutari \\
  Purdue University, USA \\
  \texttt{gscutari@purdue.edu} \\
  \And
  Alexander Rogozin \\
  MIPT and HSE University, Russia\\
  \texttt{aleksandr.rogozin@phystech.edu}
  \And
  Alexander Gasnikov \\
  MIPT, HSE University and ISP RAS\thanks{ISP RAS Research Center for Trusted Artificial Intelligence}, Russia\\
  \texttt{gasnikov@yandex.ru} \\
}

\begin{document}
 \maketitle 
 
\begin{abstract}
We study solution methods for (strongly-)convex-(strongly)-concave Saddle-Point Problems (SPPs) over networks of two type--master/workers (thus centralized) architectures  and  mesh (thus decentralized) networks.    
The local functions at each node are assumed to be \textit{similar}, due to statistical data similarity or otherwise. We establish lower complexity bounds for a fairly general  class of algorithms solving the SPP. We show that   a given suboptimality $\epsilon>0$ is achieved over master/workers networks in $\Omega\big(\Delta\cdot  \delta/\mu\cdot \log (1/\varepsilon)\big)$ rounds of communications, where $\delta>0$ measures the degree of similarity of the local functions, $\mu$ is their strong convexity constant, and $\Delta$ is the diameter of the network. The lower communication complexity bound over mesh networks reads     $\Omega\big(1/{\sqrt{\rho}} \cdot  {\delta}/{\mu}\cdot\log (1/\varepsilon)\big)$,  
where $\rho$ is the (normalized) eigengap of the gossip matrix used for the communication between neighbouring nodes.  We then propose algorithms matching the lower bounds over either  types of networks (up to log-factors). We assess the effectiveness of the proposed algorithms on a robust regression  problem.\end{abstract}

\section{Introduction}
We study smooth (strongly-)convex-(strongly-)concave  SPPs over a network of $M$  agents: 
\begin{equation}\tag{P}
\label{data_sim_intro}
    \min_{x \in X} \max_{y \in Y} f(x,y) := \frac{1}{M} \sum\limits_{m=1}^M  f_m(x,y),
\end{equation}
where $X,Y\subseteq \mathbb{R}^d$ are convex and compact sets common to all the agents; and $f_m(x,y)$ is the loss function of agent $m$, known only to the agent. Problem (\ref{data_sim_intro})  has found a wide range of applications, including, game theory 
\cite{GT-book,facchinei2007finite}, image deconvolution problems \cite{chambolle2011first},   adversarial training \cite{Arjovsky_et_al2017, Bengio2014}, and statistical learning \cite{Abadeh_et_al_2015}--see Sec.~\ref{sec:setup}  for some motivating examples in the distributed setting.   We are particularly interested in learning problems, where each $f_m$ is the empirical risk that measures the mismatch between the model to be learned and the {\it local} dataset owned by agent $m$.

Since the functions $f_m$ can be accessed only locally and routing local data to other agents is infeasible or highly inefficient, solving (\ref{data_sim_intro})  calls for the design of distributed algorithms that alternate between a local computation procedure at each agent's side, and a round of communication among (suitably chosen) neighboring nodes. We address such a design considering explicitly    two type of computational architectures, namely: {\bf (i)  master/workers} networks--these are centralized systems suitable for parallel computing; for instance, they are the typical computational architecture arising from  federated learning applications (e.g.,  \cite{kairouz2019advances}), where data are split across multiple workers and computations are performed in parallel,   coordinated by the master node(s); and \textbf{(ii) mesh networks}--these    are distributed systems with no special topology (modeled just as undirected graphs), which  capture scenarios  wherein there is  no hierarchical structure (e.g., master nodes) and each node can  communicate only   with its intermediate neighbors. 

  \textbf{Function similarity:} Motivated in particular by machine learning applications, our   design and analysis pertain to distributed algorithms for  SPPs    (\ref{data_sim_intro}) where the local functions $f_m$'s are {\it related}--quantities such as   gradients and the second derivatives matrices of   $f_m$'s differ only by  a finite quantity  $\delta>0$; we   will term such   SPPs as $\delta$-related SPPs. For instance, this is the typical situation in the aforementioned distributed empirical risk minimization setting   \cite{Arjevani-ShamirNIPS15,pmlr-v119-hendrikx20a,DISCO}: when data are i.i.d. among machines,  the $f_m$'s reflect  statistical similarities in the data residing at different nodes, resulting in a    $\delta=\tilde{O}(1/\sqrt{n})$, where $n$ is the local sample size   ($\tilde{O}$ hides log-factors and dependence on  $d$). 
  
  While SPPs have been extensively studied in the centralized setting (e.g., \cite{ facchinei2007finite, Nemirovski2004, Korpelevich1976TheEM, Nesterov07, beznosikov2021distributed}) and more recently over mesh networks  \cite{liu2019decentralizedprox, 9304470, liu2019decentralized, 7403075, rogozin2021decentralized, beznosikov2021decentralized, beznosikov2021decentralized1}, we are not aware of any analysis or (distributed) algorithm that explicitly exploit function similarity  to boost  communication efficiency--either lower complexity bounds or upper bounds. On the other hand, recent works for sum-utility minimization problems over networks (e.g.,  \cite{Arjevani-ShamirNIPS15,DANE,reddi2016aide,yuan2019convergence,GIANT,Fan2020,DISCO,pmlr-v119-hendrikx20a,sun2019distributed,NetDane})   show   that employing some form of  statistical preconditioning in the algorithm design   provably reduces communication complexity.
Whether these improvements are possible/achievable for $\delta$-related SSPs in the form (\ref{data_sim_intro})  remains unclear.  This paper provides a positive answer to the above open problem.

\textbf{Major contributions:}   Our major results are summarized next.   \textbf{(a) Lower complexity bounds:}      Under   mild structural assumptions on the algorithmic oracle  (satisﬁed by a variety of methods), we   establish lower complexity bounds for the $\delta$-related SPP (\ref{data_sim_intro}) with $\mu$-strongly-convex-strongly -concave, $L$-smooth (twice-differentiable) local functions:  an $\varepsilon$ precision on the optimality gap over master/workers system is achieved in   $\Omega\big(\Delta\cdot  \delta/\mu\cdot \log (1/\varepsilon)\big)$ communication steps, where   $\Delta$ is the diameter of the network.    The lower complexity bound over mesh networks reads $\Omega\big(1/{\sqrt{\rho}} \cdot  {\delta}/{\mu}\cdot\log (1/\varepsilon)\big)$ rounds of communications, where $\rho$ is the (normalized) eigengap of the gossip matrix used for the communication between neighbouring nodes.  These  new lower bounds show a more favorable dependence on the optimization parameters (via $\delta/\mu$) than that of distributed oracles for SPPs ignoring  function similarity \cite{beznosikov2021distributed,rogozin2021decentralized}, whose communication complexity, e.g.,  over mesh networks reads $\Omega\big(1/{\sqrt{\rho}} \cdot  {L}/{\mu}\cdot\log (1/\varepsilon)\big)$. The latter provides a pessimistic  prediction when   $\delta/\mu \ll L/\mu$.  This is the typical situation of ill-conditioned problems, such as many learning problems where  the regularization parameter that is optimal for test predictive performance is so small   that a scaling with  $L/\mu$  is no longer practical while $\delta/\mu$  is (see, e.g., \cite{marteau2019globally,pmlr-v119-hendrikx20a}). 
   \textbf{(b) Near optimal algorithms:} We proposed algorithms for such SPPs  over master/workers and mesh networks that match the lower bounds up to logarithmic factors. They are provably faster than existing solution methods for $\mu$-strongly-convex-strongly-concave, $L$-smooth SPPs, which do not exploit function similarity.  Preliminary numerical results on distributed robust logistic regression   support our theoretical findings.  \vspace{-0.4cm}
\subsection{Related works}\label{sec:Related_works}\vspace{-0.2cm}
 \textbf{Methods for SPPs ignoring function similarity:} (Strongly)-convex-(strongly)-concave SPPs have been extensively studied in  the optimization  literature   and as special instances of (strongly) monotone Variational Inequalities (VI) \cite{facchinei2007finite,juditsky2008solving}.  Several algorithms   are available in the {\bf centralized} setting, some directly imported from the VI literature; representative examples include:  the mirror-proximal algorithm \cite{Nemirovski2004}, Extragradient method \cite{Korpelevich1976TheEM} and the scheme in \cite{Nesterov07}--they are readily implementable on master/workers architectures as well.  For SPPs with $\mu$-strongly-convex-strongly-concave, $L$-smooth loss, all these schemes achieve  iteration complexity of $O\big(L/\mu \cdot\log(1/\varepsilon)\big)$, which has been shown to be optimal for first-order methods solving such a class of SPPs 
 \cite{zhang2019lower,ouyang2019lower}. 
 Lower bounds and optimal algorithms in the {\bf distributed} setting for SPPs without similarity have been studied in \cite{beznosikov2021distributed}.
 
Note that none of the above lower (and upper) complexity bounds or   (centralized or distributed) algorithmic designs capture  function similarity.  As a consequence,   convergence rates   certified in the aforementioned works, when applicable to $\delta$-related SPPs in the form (\ref{data_sim_intro}), provide quite pessimistic predictions, in the setting   $1+\delta/\mu \ll L/\mu$.

 \textbf{Methods for sum-utility minimization  exploiting function  similarity:}
Several works exploited the idea of statistical preconditioning to provably improve  communication complexity of solution methods for  the minimization of the sum of $\delta$-related, $\mu$-strongly convex and $L$-smooth functions    over  master/workers networks. Lower complexity bounds   are established in \cite{Arjevani-ShamirNIPS15}, and read $\Omega\big(\sqrt{\delta/\mu}\log(1/\varepsilon)\big)$, which contrasts with  $O\big(\sqrt{L/\mu}\log(1/\varepsilon)\big)$ achievable by first-order (Nesterov) accelerated methods \cite{nesterov2018lectures}, certifying thus   faster rates   whenever  $\delta/\mu<L/\mu$. Solutions methods exploiting function similarity  are   mirror proximal-like schemes, and  include \cite{DANE,reddi2016aide,yuan2019convergence} (for quadratic losses), \cite{DISCO} (for self-concordant losses),   \cite{GIANT}, and  \cite{Fan2020} (for composite optimization), with  \cite{pmlr-v119-hendrikx20a} employing acceleration.   None of these methods are implementable over  mesh networks, because they   rely on  a centralized (master) node. To our knowledge, Network-DANE \cite{NetDane} and SONATA \cite{sun2019distributed} are the only two methods   that leverage statistical similarity to enhance convergence of distributed methods over mesh networks; \cite{NetDane} studies  strongly convex quadratic losses while \cite{sun2019distributed} considers general objectives, achieving a communication complexity of $\widetilde{{O}}((1/\sqrt{\rho})\cdot\delta/\mu\cdot \log(1/\varepsilon))$, where $\widetilde{{O}}$ hides logarithmic factors. None of the methods above however are applicable to the $\delta$-related SPP (\ref{data_sim_intro}).

\subsection{Notation}
Given a positive integer $M$, we define $[M]=\{1,\ldots, M\}$.  We use $\langle x,y \rangle := \sum_{i=1}^d x_i y_i$ to denote standard inner product of $x,y\in\R^d$. It induces $\ell_2$-norm in $\R^d$ in the following way $\|x\| := \sqrt{\langle x, x \rangle}$. We also introduce $\text{proj}_{\mathcal{Z}}(z) = \min_{u \in \mathcal{Z}}\| u - z\|$ -- the Euclidean projection onto $\mathcal{Z}$. We order the eigenvalues of any symmetrix matrix $A\in \mathcal{R}^{m\times m}$ in nonincreasing fashion, i.e., $\lambda_{\max}(A)=\lambda_1(A)\geq \ldots \geq  \lambda_m(A)=\lambda_{\min}(A)$, with   $\lambda_{\max}(\cdot)$ [resp. $\lambda_{\min}(\cdot)$] denoting the largest (resp. smallest) eigenvalue.

\section{Setup and Background}\label{sec:setup}

\textbf{Problem setting:} We begin introducing the main assumptions underlying Problem~\eqref{data_sim_intro} and some useful notation. 
 
 Let us stack the $x$- and $y$-variables in the tuple $z = (x,y)$; accordingly,   define $\mathcal{Z} = \mathcal{X} \times \mathcal{Y}$ and    
  the  vector-functions $F_m, F:\mathcal{Z}\to \mathcal{R}^{2d}$:
\begin{equation}
\label{opSP}
    F_m(z) :=  \begin{pmatrix}
\nabla_x f_m(x,y)\\
-\nabla_y f_m(x,y)
\end{pmatrix},\quad \text{and} \quad F(z):= \frac{1}{M} \sum\limits_{m=1}^M F_m(z).
\end{equation}
 
 The following conditions are standard for strongly convex-strongly concave SPPs.
 \begin{assumption}\label{ass_problem}  Given \eqref{data_sim_intro}, the following hold:
 \begin{itemize}
 \item[(i)] $\emptyset\neq \mathcal{Z}$ is a convex set; 
     \item[(ii)] Each $f_m:\mathcal{R}^{2d}\to \mathcal{R}$ is twice differentiable on (an open set containing) $\mathcal{Z}$, with $L$-Lipschitz gradient: 
   $ \|F_m(z_1) - F_m(z_2)\| \leq L\|z_1-z_2\|$, for all   $z_1, z_2\in \mathcal{Z}$;
     \item[(iii)]  $f(z)$ is  $\mu$-strongly convex-strongly concave on $\mathcal{Z}$, i.e., 
      $  \langle F(z_1) - F(z_2), z_1 - z_2 \rangle \geq \mu\|z_1-z_2\|^2$, for all  $z_1, z_2\in \mathcal{Z}$;
      \item[(iv)]  Each $f_m(z)$ is  convex-concave on $\mathcal{Z}$, i.e. $0$-strongly convex-strongly concave.
 \end{itemize}
 
 \end{assumption}
 
 We are interested in finding the solution $z^* = (x^*, y^*)$ of Problem~\eqref{data_sim_intro}  under function similarity.  
 \begin{assumption}[$\delta$-related $f_m$'s] \label{ass_similarity}The local functions are $\delta$-related: for all $(x,y)\in \mathcal{Z}$, 
\begin{align*} 
    \| \nabla^2_{xx} f_{m} (x,y) -  \nabla^2_{xx} f (x,y) \| \leq \delta, \notag\\ \| \nabla^2_{xy} f_{m} (x,y) -  \nabla^2_{xy} f (x,y) \| \leq \delta, \\ \| \nabla^2_{yy} f_{m} (x,y) -  \nabla^2_{yy} f (x,y) \| \leq \delta \notag.
\end{align*}\end{assumption}
The interesting case is when $1+\delta /\mu \ll L/\mu$. When the $f_m$'s are empirical loss functions over local data sets of size $n$,  under standard assumptions on data distributions and learning model (e.g., \cite{DISCO,pmlr-v119-hendrikx20a}),  
$\delta=\tilde{O}(1/\sqrt{n})$ with high probability  ($\tilde{O}$ hides log-factors and dependence on  $d$)--some motivating examples falling in this category are discussed in Sec.~\ref{sec_examples} below.   While  such examples represent important applications, we point out that  our (lower and upper) complexity bounds are valid in all scenarios wherein Assumption~\ref{ass_similarity}  holds, not necessarily due to  statistical arguments. 

\textbf{Network setting:} 
The communication network  is modeled as a fixed, connected, undirected graph,   $\mathcal{G} \triangleq (\mathcal{V}, \mathcal{E})$, where $\mathcal{V} \triangleq \{1,\ldots, M\}$ denotes  the vertex    set--the set of agents--while $\mathcal{E} \triangleq \{(i,j) \, |\, i,j \in \mathcal{V} \}$ represents the set of edges--the communication links;  $(i,j) \in \mathcal{E}$ iff there exists a communication link between agent $i$ and $j$.  We denote by $\Delta$ the diameter of the graph. 
When it comes to distributed algorithms over mesh networks, we leverage neighbouring communications among adjoining nodes.  Communications of $d$-dimensional vectors  will be modeled as a  matrix multiplication by a matrix $W$ (a.k.a. gossip matrix). The following assumptions on $W$ are standard to establish convergence of   distributed algorithms over mesh networks. 
\begin{assumption}\label{Ass_gossip_matrix}
    The matrix $W\in \mathcal{R}^{M\times M}$   satisfies the following: 
    \textbf{(a)} It is compliant with $\mathcal{G},$ that is, (i) $w_{ii}>0, \forall i\in [M];$ 
(ii) $w_{i j}>0,$ if $\{j,i\}\in \mathcal{E};$ and (iii)  $w_{i j}=0$ otherwise;
 \textbf{(b)} It is symmetric and   stochastic, that is,    $W1=1$ (and thus also  $1^{\top}W=1^{\top}$).
\end{assumption}

 Notice that  a direct consequence of   Assumption~\ref{Ass_gossip_matrix} (along with the fact that $\mathcal{G}$ is connected) is that
\begin{equation}\label{connectivity number}
 \rho\triangleq 1-\max\{\lambda_2(W), \lvert\lambda_{\min}(W)\rvert\}<1, 
 \end{equation} where $\rho$ is the eigengap between the first and second largest (magnitude) eigenvalue of $W$.    Roughly speaking, $\rho$  measures  how fast the network mixes information (the larger, the faster).

\subsection{Motivating examples}
\label{sec_examples}

Several problems of interest can be cast in the SPP (\ref{data_sim_intro}), for which function similarity arises naturally, some are  briefly discussed next.

\textbf{Robust Regression}: Consider  the robust instance of the linear regression problem in its Lagrangian form:  
\begin{align}\label{eq:regression}
    \min_{w} \max_{r} \frac{1}{2N} \sum\limits_{i=1}^N (w^T (x_i + r) - y_i)^2 + \frac{\lambda}{2} \| w\|^2 - \frac{\beta}{2} \|r \|^{2},
\end{align}
where $w$ are the weights of the model, $\{(x_i, y_i)\}_{i=1}^N$ are pairs of the training data, and $r$ models the noise, and   $\lambda$ and $\beta$ are the regularization parameters. Let $n$ be the local sample size (thus $N=nm$).  The typical regularization parameter that is optimal for test predictive performance is  $\lambda=O(1/\sqrt{N})$. Assuming $\beta$ of the same order of $\lambda$ and  invoking function similarity $\delta=O(1/\sqrt{n})$    \cite{marteau2019globally,pmlr-v119-hendrikx20a} yield a condition number of the problem $\kappa=\mathcal{O}(\sqrt{m\cdot n})$ while $\delta/\mu=\mathcal{O}(\sqrt{m})$. This implies that  first order methods applied to \eqref{eq:regression} will slowdown as the local sample size $n$ grows. Rate scaling with $\delta/\mu$ would be instead  independent on the local sample size.

\textbf{Adversarial robustness of neural networks:} Recent works have demonstrated that deep neural networks are vulnerable to adversarial examples—inputs that are almost indistinguishable from natural data and yet classified incorrectly by the network \cite{szegedy2013intriguing,goodfellow2014explaining}. To improve   resistance   to a variety  of adversarial inputs, a  widely studied approach is leveraging robust optimization and formulate the training as saddle-point problem \cite{madry2017towards,nouiehed2019solving}:
\begin{align*}
    \min_{w} \max_{r} \frac{1}{N} \sum\limits_{i=1}^N l(f(w, x_i + r, y_i)^2 + \frac{\lambda}{2} \| w\|^2 - \frac{\beta}{2} \|r \|^{2},
\end{align*}
where $w$ are the weights of the model, $\{(x_i, y_i)\}_{i=1}^N$ are pairs of the training data, $r$ is the so-called adversarial noise, which models a  perturbation in the data, and   $\lambda$ and $\beta$ are the regularizers.

\textbf{Other optimization problems:} Other instances of the SPP are the (online) transport or Wasserstein Barycenter (WB) problems, see \cite{jambulapati2019direct,dvinskikh2021improved}. This representation comes from the dual view of transportation polytope. 
b) Another example  is Lagrangian based optimization problems. For instance, consider the minimization of the sum of loss functions, each one associated to one agent, subject to some (common) constraints. The problem can be equivalently    rewritten as a saddle-point problem using Lagrangian multipliers. It is easy to check that if the agents' functions are $\delta$-related, then the resulting  saddle-point problem is also so.

\section{Lower Complexity Bounds}\label{sec:lower_bound} 
In this section we establish lower complexity bounds for centralized (i.e., master/workers-based) and distributed (gossip-based) algorithms. We begin introducing  the back-box procedure describing the class of algorithms these lower bounds pertain to.      
\subsection{Optimization/communication oracle} \label{oracle} 
Our   procedure models a fairly general class of (centralized and distributed) algorithms over graphs, whereby nodes perform local computation  and communication  tasks. Computations at each node are  based on linear operations involving current or past iterates, gradients, and vector products with local Hessians and their inverses, as well as solving local optimization problems involving such quantities. During communications, the nodes can share (compatibly with the graph topology)  any of the vectors they have computed up until that time. The black-box procedure can be formally describe as follows.  

\begin{definition}[Oracle] \label{proc}
    Each agent $m$ has its own local memories $\mathcal{M}^x_{m}$ and $\mathcal{M}^y_{m}$ for the $x$- and $y$-variables, respectively--with  initialization $\mathcal{M}_{m}^x = \mathcal{M}_{m}^y= \{0\}$.    $\mathcal{M}_{m}^x$ and $\mathcal{M}_{m}^x$ are updated as follows.\\
   
    $\bullet$ \textbf{Local computation:} Between communication rounds, each agent $m$ computes and adds to its $\mathcal{M}^x_{m}$ and $\mathcal{M}^y_{m}$ a finite number of points $x,y$, each satisfying 
    \begin{equation}\begin{aligned}\label{eq:oracle-opt-step}
        \alpha x + \beta \nabla_x f_m(x,y) \in \text{span} \big\{&x'~,~\nabla_x f_m(x',y'),\\ 
        &(\nabla^2_{xx} f_m (x'',y'') + D)x'~,~ (\nabla^2_{xx} f_m (x'',y'') + D)\nabla_x f_m(x',y') \\ &(\nabla^2_{xx} f_m (x'',y'') + D)^{-1}x'~,~ (\nabla^2_{xx} f_m (x'',y'') + D)^{-1}\nabla_xf_m(x',y'), \\
        &(\nabla^2_{xy} f_m (x'',y''))y' ~,~  (\nabla^2_{xy} f_m (x'',y''))\nabla_y f_m(x',y') \big\},\\
        \theta y - \varphi \nabla_y f_m(x,y) \in \text{span} \big\{&y'~,~\nabla_y f_m(x',y'),\\
        &(\nabla^2_{yy} f_m (x'',y'') + D)y'~,~ (\nabla^2_{yy} f_m (x'',y'') + D)\nabla_y f_m(x',y') \\ &(\nabla^2_{yy} f_m (x'',y'') + D)^{-1}y'~,~ (\nabla^2_{yy} f_m (x'',y'') + D)^{-1}\nabla_y f_m(x',y'), \\
        &(\nabla^2_{xy} f_m (x'',y''))^T x' ~,~  (\nabla^2_{xy} f_m (x'',y''))^T \nabla_x f_m(x',y') \big\}, 
    \end{aligned}\end{equation}
    for given $x', x'' \in \mathcal{M}^x_{m}$ and  $y', y'' \in \mathcal{M}^y_{m}$;  some $\alpha,\beta,\theta,\varphi \geq 0$ such that  $\alpha + \beta >0$ and  $\theta + \varphi >0$; and   $D$ is some diagonal matrix (such that all the inverse matrices   exist).

    $\bullet$ \textbf{Communication:} Based upon communication rounds  among neighbouring nodes,   $\mathcal{M}^x_{m}$ and $\mathcal{M}^y_{m}$ are updated according to
    \begin{equation}\label{eq:oracle-comm}
        \mathcal{M}^x_{m} := \text{span}\left\{\bigcup_{(i,m) \in \mathcal{E}} \mathcal{M}^x_{i} \right\}, \quad 
        \mathcal{M}^{y}_{m} := \text{span}\left\{\bigcup_{(i,m) \in \mathcal{E}} \mathcal{M}^y_{i} \right\}.
    \end{equation}

    $\bullet$ \textbf{Output:} 
    The final global output is calculated as: 
    \begin{align*}
        x^{K} \in \text{span}\left\{\bigcup_{m=1}^M \mathcal{M}^x_{m} \right\},~~y^{K} \in \text{span}\left\{\bigcup_{m=1}^M \mathcal{M}^y_{m} \right\}.
    \end{align*}
\end{definition} 

The above oracle captures a gamut of existing centralized and distributed algorithms.   For instance,  local computations  model  either inexact local solutions--e.g., based on single/multiple steps of  gradient or Newton-like updates, which corresponds to setting $\alpha=\theta=1$ and $\beta=\varphi=0$--or exact solutions of agents' subproblems (via some subroutine algorithm), corresponding to $\alpha=\theta=0$ and $\beta=\varphi=1$. Multiple     rounds of computations (resp. communications)  can be performed between communication rounds  (resp. computation tasks). 
Notice that the proposed oracle builds on   \cite{scaman2017optimal,Arjevani-ShamirNIPS15}  for minimization problems over networks--the former modeling only gradient updates  and the latter considering only centralized optimization (master/workers systems). 

\subsection{Lower complexity bounds}
We are in the position to state our main results on lower communication complexity--Theorem~\ref{th-LB-centralized}  pertains to algorithms over master/workers systems while Theorem~\ref{th-LB-distributed} deals with mesh networks.
\begin{theorem}\label{th-LB-centralized}
For any $L, \mu, \delta >0$ and connected graph $\mathcal G$ with diameter  $\Delta>0$, there exist a   SPP in the form (\ref{data_sim_intro})  (satisfying Assumption \ref{ass_problem}) with $\mathcal{Z}=\mathcal{R}^{2d}$ (where $d$ is sufficiently
large), $x^* \neq 0$, $y^* \neq 0$,  and  local functions $f_m$ being  $L$-smooth, $\mu$-strongly-convex-strongly-concave, $\delta$-related (Assumption \ref{ass_similarity}) such that any centralized algorithm satisfying Definition~\ref{proc} produces    the following estimate on the global output $z^K=(x^K,y^K)$ after $K$ communication rounds:
\begin{equation*}
    \|z^{K} - z^*\|^2 = \Omega\left(\exp\left( - \frac{K}{\Delta } \cdot \frac{1}{\frac{1}{8}\sqrt{1 + \left(\frac{\delta}{32\mu}\right)^2} - \frac{1}{8}} \right) \| y^*\|^2\right).
\end{equation*}
\end{theorem}

\begin{corollary}
In the setting of Theorem~\ref{th-LB-centralized}, 
the number of communication rounds required to obtain a $\varepsilon$-solution is lower bounded by
\begin{equation}
\label{cntr_b}
    \Omega\left( \Delta\left(1 + \frac{\delta}{\mu}\right) \cdot  \log \left(\frac{\| y^*\|^2}{\varepsilon}\right)\right).
\end{equation}
\end{corollary}

\begin{theorem}\label{th-LB-distributed}
For any $L, \mu, \delta >0$ and $\rho \in (0;1]$, there exist a  SPP in the form (\ref{data_sim_intro})  (satisfying Assumption \ref{ass_problem}) with $\mathcal{Z}=\mathcal{R}^{2d}$(where $d$ is sufficiently
large), $x^* \neq 0$, $y^* \neq 0$, and  local functions $f_m$ being  $L$-smooth, $\mu$-strongly-convex-strongly-concave, $\delta$-related (Assumption \ref{ass_similarity}), and a gossip matrix $W$ over the connected graph $\mathcal{G}$,  satisfying Assumption \ref{Ass_gossip_matrix} and  with eigengap  $\rho$,      such that any decentralized algorithm  satisfying Definition~\ref{proc} and using the gossip matrix $W$ in the communication steps \eqref{eq:oracle-comm}  produces    the following estimate on the global output $z^K=(x^K,y^K)$ after $K$ communication rounds:
\begin{equation*}
    \|z^{K} - z^*\|^2 = \Omega\left(\exp\left(- K \sqrt{\rho} \cdot \frac{1}{\frac{1}{20}\sqrt{1 + \left(\frac{\delta}{32\mu}\right)^2} - \frac{1}{20}} \right)  \| y^*\|^2\right).
\end{equation*}
\end{theorem}

\begin{corollary}
In the setting of Theorem~\ref{th-LB-distributed}, the number of communication rounds required to obtain a $\varepsilon$-solution is lower bounded by
\begin{equation}
\label{decntr_b}
    \Omega\left( \frac{1}{\sqrt{\rho}}\left(1 + \frac{\delta}{\mu}\right) \cdot  \log \left(\frac{\| y^*\|^2}{\varepsilon}\right)\right).
\end{equation}
\end{corollary}

  These lower  complexity bounds   show an expected   dependence   on the  optimization parameters and network quantities.  Specifically, the number of communications scale proportionally  to  $\delta/\mu$--this generalizes existing lower bounds \cite{beznosikov2021distributed} that do not account for such similarity, resulting instead in the more pessimistic dependence on $L/\mu$--typically $\delta\leq L$. The network impact is captured by  the diameter  $\Delta$ of the network for master/workers architectures--$\Delta$ communications  steps are required in the worst case to transmit a message between two nodes--and the eigengap $\rho$ of the matrix $W$, when arbitrary graph typologies are consider; $1/\sqrt{\rho}$  can be bounded as $\mathcal{O}(T)$, where $T$ is the largest hitting time of the Markov chain with  probability transition matrix  $W$  \cite{Olshevsky}.  For instance, for fully connected networks $\Delta=1/\sqrt{\rho}=1$ while for star networks $\Delta=1$ and $1/\sqrt{\rho}=\sqrt{M}$.  For general graphs, $1/\sqrt{\rho}$ can be larger than $\Delta$, see 
   \cite{Nedic_Olshevsky_Rabbat2018} for more details.  
To certify the tightness of the derived lower bounds, the next section designs    algorithms that reach such bounds.  

\section{Optimal algorithms}
\subsection{Centralized case (master/workers systems)}

Our first optimal algorithm is for SPPs over master/workers architectures or more generally networked systems where    a spanning tree (with the root as master node) is preliminary set; it is formally described in Algorithm~\ref{alg:cent}. We assumed w.l.o.g. that the  master node owns function $f_1$. 

Some insights on the genesis of this method are discussed next. 

$\bullet$ Consider for a moment the minimization problem $\min_{x \in X} f(x) := \frac{1}{M} \sum_{m=1}^M  f_m(x)$, under Assumption 2. Following \cite{DANE} we can solve it invoking the mirror descent algorithm, which reads  \begin{equation}
    \label{step_min}
    x^{k+1} = \arg\min_x \left[\langle \eta \nabla f(x^k), x \rangle + D_{\phi}(x, x^k)\right],
\end{equation}
where $D_{\phi}(x, y) = \phi(x) - \phi(y) - \langle \nabla \phi(y), x-y \rangle$ is the Bregman divergence, with function $\phi(x) = f_1(x) + \frac{\delta}{2}\|x\|^2$. It is shown that we can take stepsize $\eta = 1$ (\cite{zhang2018communication, pmlr-v119-hendrikx20a}). Therefore,   \eqref{step_min} can be rewritten as  
\begin{equation}
    \label{step_min_new}
    x^{k+1} = \arg\min_x \left[\frac{1}{\delta} f_1(x) + \frac{1}{2} \left\|x - x^k + \frac{1}{\delta} (\nabla f(x^k) - \nabla f_1(x^k)) \right\|^2\right].
\end{equation}
Noting that in Algorithm 1 $\gamma \sim \frac{1}{\delta}$ (see Appendix B.1), one infers the connection  between  \eqref{step_min_new} and the updates in  lines 3~(i) and 3~(ii). The extra step as in  line 3~(iii) is due to the fact that Algorithm 1 solves a SPP (and not a classical minimization as postulated above):   gradient descent-like methods as  \eqref{step_min} are not  optimal   for SPPs; in fact, they might  diverge when applied to  general convex-concave SPPs. Out approach is then to employ  Forward-Backward-Forward algorithms \cite{tseng2000modified} or the Extragradient \cite{Korpelevich1976TheEM} method, which leads  to the step in  line 3~(iii).

$\bullet$ Another interpretation of the proposed algorithm comes from  looking at Problem \eqref{data_sim_intro} as a composite minimization problem, with objective function $h_1(x,y) + h_2(x,y)$, with $h_1(x,y) = f_1(x,y)$ and $h_2(x,y) = \frac{1}{M} \sum_{m=1}^M (f_m(x,y) - f_1(x,y))$. The  first function $h_1$ is $L$-smooth and convex-concave while  $h_2$ is $\delta$-smooth and, in general, non-convex-non-concave. Such type of problems can be solved invoking sliding techniques \cite{lan2016gradient,rogozin2021decentralized}. 

\begin{algorithm}[th]
	\caption{\bf (Star Min-Max Data Similarity Algorithm)}
	\label{alg:cent}
	\hspace*{\algorithmicindent} {\bf Parameters:} stepsize $\gamma$, accuracy $e$;\\
	\hspace*{\algorithmicindent} {\bf Initialization:} Choose  $z^0=(x^0,y^0)\in \mathcal{Z}$, $z^0_m = z^0$, for all $m\in [M]$;
	\begin{algorithmic}[1]
		\For{$k=0,1,2,\ldots$}
		\State  {Each worker}  $m$ computes $F_m (z^k)$ and  {sends} it to the  master; \label{line2}
		\State  {The master node}:  \begin{itemize}
		    \item[(i)] computes $v^k = z^k - \gamma \cdot \left(F(z^k) - F_1(z^k)\right)$;
		    \item[(ii)]   finds $u^k$, s.t. $ \| u^k - \hat{u}^k \|^2 \leq e$, where $\hat{u}^k$ is the solution of: \label{line4}
		\begin{equation}
		\label{inter}
			\min_{u_x \in \mathcal{X}} \max_{u_y \in \mathcal{Y}}  \left[\gamma f_1(u_x,u_y) + \frac{1}{2} \|u_x- v^k_x \|^2 - \frac{1}{2} \|u_y- v^k_y \|^2\right];
		\end{equation}
		\item[(iii)] updates $z^{k+1} = \text{proj}_{\mathcal{Z}} \left[u^k + \gamma \cdot (F(z^k) - F_1(z^k) - F(u^k) + F_1(u^k))\right]$ and broadcasts $z^{k+1}$ to the workers
		\end{itemize} 
		\EndFor
	\end{algorithmic}
\end{algorithm}

It is not difficult to check that    Algorithm~\ref{alg:cent} is an instance of the oracle introduced in Definition~\ref{proc}. It accommodates either exact solutions of the strongly convex  subproblems~\eqref{inter} (corresponding to $e=0$) or inexact ones (up to tolerance $e>0$)--the latter can be computed, e.g., using Extragradient method \cite{juditsky2008solving}, which is optimal in this case.

The communication complexity of the method is proved in the next theorem, which certifies that the proposed algorithm is optimal, i.e.,  achieves the lower bound (\ref{cntr_b}) on the number of required  communications--we refer to Appendix~\ref{app:op_c} in the supplementary material for a detailed description of the algorithmic tuning as well as a study of the computational complexity when Extragradient method is employed to solve  subproblems~\eqref{inter} (up to a suitably chosen tolerance). 

\begin{theorem} \label{theorem2}
Consider Problem~\eqref{data_sim_intro} under Assumptions~\ref{ass_problem}-\ref{ass_similarity} over a connected graph $\mathcal{G}$ with a master node. Let  $\{z^k\}$ be the sequence generated by Algorithm~\ref{alg:cent} with tuning as described in  Appendix~\ref{app:op_c} (cf.  the supplementary material). Then, given $\varepsilon>0$, the number of communication rounds for 
$\|z^k - z^*\|^2 \leq \varepsilon$   is  $\cO \big(\big(1 + {\delta}/{\mu}\big)\log({1}/{\varepsilon})\big)$. 
\end{theorem}

\subsection{Distributed case (mesh networks)}

We consider now mesh networks. Because of the lack of a master node, each agent $m$ now owns local estimates $u_m$ and $v_m$ of the common variables  $u$ and $v$, respectively, which are iteratively updated.  At each iteration,   a node is selected uniformly at random, which plays the role of the master node, performing thus the  update of its own local variables, followed by some rounds of communications via  accelerated  (inexact) gossip protocols \cite{liu2011accelerated,ye2020multi}--the latter being instrumental to propagate the updates of the $u,v$-variables and gradients across the network.     The algorithm is formally introduced in Algorithm~\ref{alg:decent}, with the accelerated gossip procedure described in Algorithm~\ref{accgossip}.    
\begin{algorithm}[H]
	\caption{\bf (Distributed  Min-Max Data Similarity Algorithm)}
	\label{alg:decent}
	\hspace*{\algorithmicindent} {\bf Parameters:} stepsize $\gamma$, accuracy $e, e_0, e_1$, communication rounds $H_0$, $H_1$;\\
	\hspace*{\algorithmicindent} {\bf Initialization:} Choose  $z^0=(x^0,y^0)\in \mathcal{Z}$, $z^0_m = z^0$, for all $m\in [M]$;
	\begin{algorithmic}[1]
		\For{$k=0,1,2,\ldots$}
			\State  {\texttt{Communications}:} $\bar F^k_1, \ldots \bar F^k_M  = \text{AccGossip}(F_1 (z^k_1),\ldots F_M (z^k_M); H_0)$;\label{lin4}\smallskip 
		\State \texttt{Local computations}: Choose an index $m_k\in [M]$ uniformly at random; then node $m_k$  \label{lin2}
	 \begin{itemize}
	     \item[(i)] computes $v^k_{m_k} = z^k_{m_k} - \gamma \cdot (\bar F^k_{m_k} - F_{m_k} (z^k_{m_k}))$;
	     \item[(ii)] finds $\tilde u^k_{m_k}$, s.t. $ \| \tilde u^k_{m_k} - \hat{u}^k_{m_k} \|^2 \leq e$, where $\hat{u}^k_{m_k}$ is the solution of:\begin{equation}\label{inner}
			\min_{u_x\in \mathcal{X}} \max_{u_y\in \mathcal{Y}}  \left[\gamma f_{m_k}(u_x,u_y) + \frac{1}{2} \|u_x- v^k_{x,{m_k}} \|^2 - \frac{1}{2} \|u_y- v^k_{y,{m_k}} \|^2\right];
		\end{equation}
	 \end{itemize}
	 \State \texttt{Communications:} Run accelerated gossip to propagate $\tilde u^k_{m_k}$  and update  gradient variables: 
	\begin{align*}
	   &  u^k_1, \ldots u^k_M = M \cdot \text{AccGossip}( 0, \ldots, 0, \tilde u^k_{m_k}, 0 \ldots, 0; H_1),\\ 
	 &    \bar F^{k+1/2}_1,\ldots\bar F^{k+1/2}_M = \text{AccGossip}(F_1 (u^k_1), \ldots F_M (u^k_M); H_0);
	\end{align*}  
	 \State {\texttt{Update} of $\tilde z_{m_k}$-variable:} node $m_k$ performs 
		$$\tilde z^{k+1}_{m_k} =  \tilde u^k_{m_k} + \gamma \cdot (\bar F^k_{m_k} - F_{m_k} (z^k_{m_k}) - \bar F^{k+1/2}_{m_k} + F_{m_k} (\tilde u^k_{m_k}));$$
	\State \texttt{Communications:} Run accelerated gossip to propagate $\tilde z^{k+1}_{m_k}$: $$\hat z^{k+1}_1,\ldots \hat z^{k+1}_M= M \cdot \text{AccGossip}(0 \ldots, 0, \tilde z^{k+1}_{m_k}, 0 \ldots, 0; H_1);$$  
		\State  {Each worker} update $z^{k+1}_i = \text{proj}_{\mathcal{Z}} \left[\hat z^{k+1}_i\right];$
		\EndFor
	\end{algorithmic}
\end{algorithm}
\begin{algorithm} [H]
	\caption{{\bf (AccGossip)}}
	\label{accgossip}
	\hspace*{\algorithmicindent} {\bf Input:}   $z_1, ..., z_M\in \mathcal{R}^{2d}$, and  $H>0$ (communication rounds); \\
	\hspace*{\algorithmicindent} {\bf Initialization:} Construct matrix $Z$ with rows $z^T_1, \ldots, z^T_M$; Set\\
\hspace*{2.7cm}	  $Z^{-1}=Z$, \quad  $Z^0 = Z$, \quad  and \quad $\eta = \frac{1 - \sqrt{1 - \lambda_2^2( W)}}{1 + \sqrt{1 - \lambda_2^2(W)}}\vspace*{-0.3cm}$. 
	\begin{algorithmic}[1]
 
\For {$t=0,1, 2, \ldots, H$ } \smallskip

\State $Z^{t+1} = (1 + \eta) W Z^{t} - \eta Z^{t-1}$,
\EndFor
	\end{algorithmic}
	\hspace*{\algorithmicindent} {\bf Output:} Rows of  $Z^{H+1}$
\end{algorithm}

 Convergence of the method is established in Theorem~\ref{theorem3}  below--we refer to Appendix~\ref{app:op_dec} in the supplementary material for a detailed description of the algorithmic tuning  [choice of the stepsize $\gamma$, precision $e$, numbers of communications rounds $H_0, H_1$, and algorithm to solve  \eqref{inter}].

\begin{theorem}  \label{theorem3} 
Consider Problem~\eqref{data_sim_intro} under Assumptions~\ref{ass_problem}-\ref{ass_similarity} over a connected graph $\mathcal{G}$. Let  $\{(z^k_m)_{m\in [M]}\}$ be the sequence generated by Algorithm~\ref{alg:decent} with tuning as described in  Appendix~~\ref{app:op_dec}  (cf.  the supplementary material) and gossip matrix $W$ satisfying Assumption~\ref{Ass_gossip_matrix}. Then, given $\varepsilon>0$, the number of communication rounds for $\|\bar z^k - z^*\|^2 \leq \varepsilon$   reads
$\mathcal{\tilde O} \left( {1}/{\sqrt{\rho}} \cdot  \big(1 + {\delta}/{\mu}\big)\log^2\frac{1}{\varepsilon}\right)$, where $\bar z^k = \frac{1}{M} \sum_{m=1}^M z^k_m$.
\end{theorem}
 While the algorithm achieves the lower bound (\ref{decntr_b}), up to log-factors (which now however depends on $\varepsilon$ as well), there is room for improvements. In fact,   selecting  only one agent at time performing the updates   does not fully exploit the potential computational speedup offered by the networking setting. Also, the use of   gossip protocols to propagate the updates of a single agent across the entire network seems to be not quite efficient.  Designing alternative distributed algorithms overcoming these limitation is a challenging open problem.  

\section{Numerical Results}\label{sec:numerical_results}

We simulate the Robust Linear Regression problem which is defined as
\begin{align}\label{eq:rob-regr}
    \min_{w} \max_{\norm{r}\leq R_r} \frac{1}{2N} \sum\limits_{i=1}^N (w^T (x_i + r) - y_i)^2 + \frac{\lambda}{2} \| w\|^2 - \frac{\beta}{2} \|r \|^2.
\end{align}
where $w$ are the model weights, $\{x_i, y_i\}_{i=1}^N$ is the training dataset, and $r$ is the artificially added noise; we use $\ell_2$-regularization  on both $w$ and $r$. 
We solve the problem  over a master/workers topology; we consider a network with 25 workers. We test Algorithm~1 wherein the subproblems \eqref{inter} at the master node  are solved  with high accuracy  using Extragradient method. A description of the tuning of the algorithm parameters can be found in Appendix~\ref{app:numerical_results}. The algorithms are implemented in Python 3.7\footnote{Source code: \href{https://github.com/alexrogozin12/data_sim_sp}{https://github.com/alexrogozin12/data\_sim\_sp}}.
\begin{figure}[ht]
    \includegraphics[width=1\textwidth]{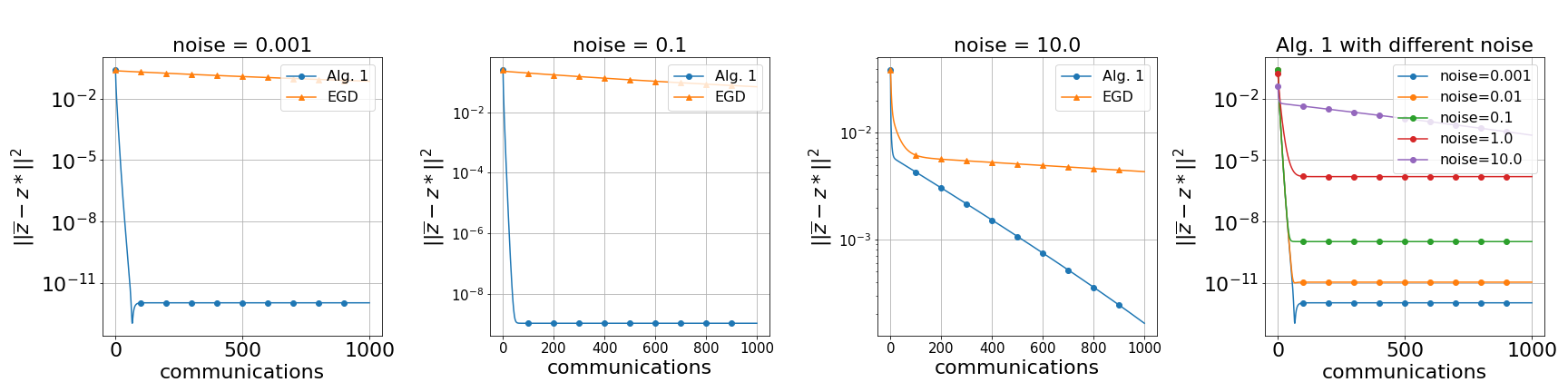}
    \caption{Centralized case, simulated data, 25 workers, ambient dimension = $40$}
    \label{fig:sim_2}
\end{figure}

\begin{figure}[H]
    
    \begin{subfigure}{0.48\textwidth}
        \centering
        \includegraphics[width=1.0\textwidth]{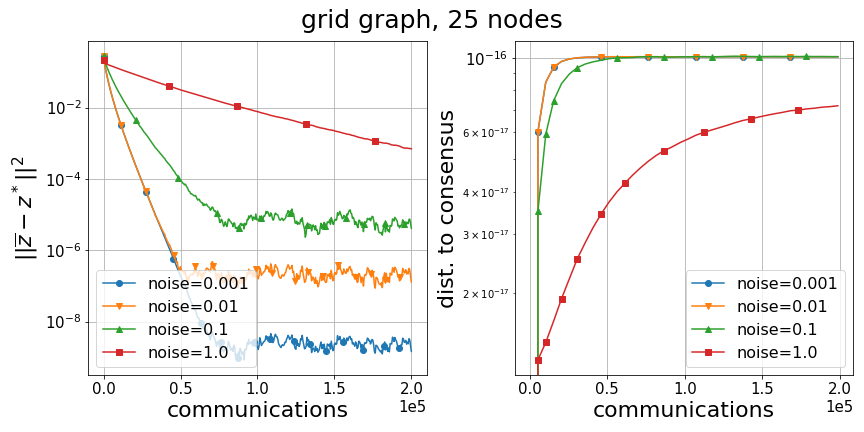}
        \label{subfig:grid}
    \end{subfigure}
    \hspace{0.3cm}
    \begin{subfigure}{0.48\textwidth}
        \centering
        \includegraphics[width=1.0\textwidth]{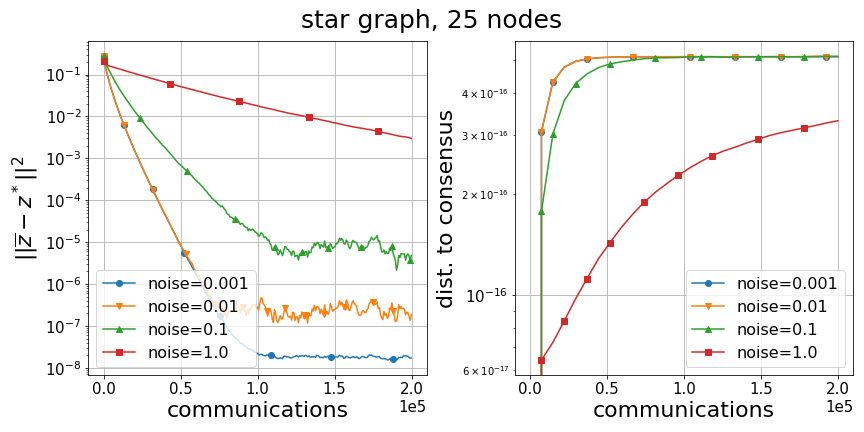}
        \label{subfig:star}
    \end{subfigure}
    \caption{Decentralized case, Alg. 2 with different noise}
    \label{fig:decentralized_sliding_different_noise}
\end{figure}

\begin{wrapfigure}[8]{r}{0.35\textwidth}
    \vspace{-0.3cm}
    \includegraphics[width=0.35\textwidth]{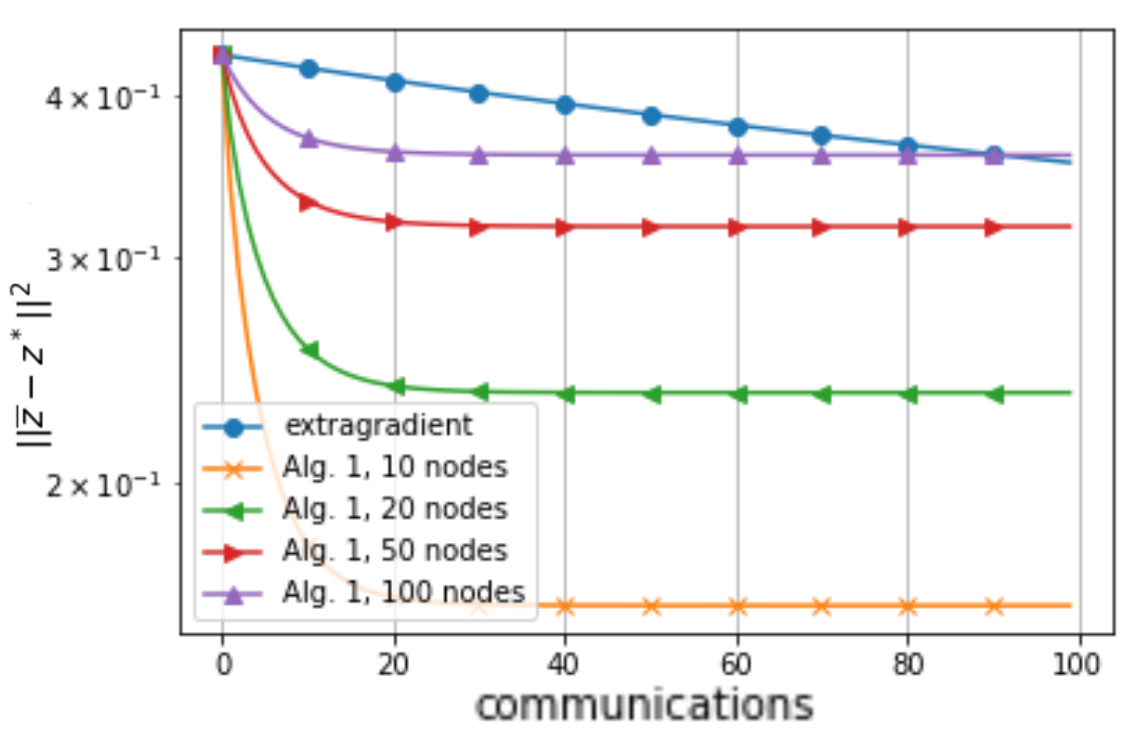}
\caption{Centralized case, a9a dataset}\label{fig:a9a_centralized}
\end{wrapfigure}

Our first experiment uses  synthetic data, which allows us to control the factor $\delta$, measuring statistical similarity of functions over different nodes. Specifically, we assume all local  datasets of size  {$n=100$}. The data set $\{\hat x_i, \hat y_i\}_{i=1}^{ n}$ at the master node is generated  randomly, with  each entry of $\hat x_i$ and  $\hat y_i$, $i=1, \ldots, n$ drawn from  the Standard Gaussian distribution. The datasets at the workers' sides,  $i = 2, \ldots, M$, are obtained perturbing $\{\hat x_i, \hat y_i\}_{i=1}^{ n}$ by  random noise $\xi_i$ with controlled variance.


Figure~\ref{fig:sim_2} compares the performance of Algorithm~\ref{alg:cent} and the Centralized Extragradient method \cite{beznosikov2021distributed} applied to Problem (\ref{eq:rob-regr}), under  different level of noise added to local datasets (level of similarity), and two different problem and network dimensions -- we plot the distance of the iterates from the solution versus the number of communications. It can be seen that Algorithm~\ref{alg:cent} consistently outperforms the Extragradient method in terms of number of communications--the smaller the noise (the more similar the local functions are), the larger the gap between the two algorithm (in favor of Algorithm~\ref{alg:cent}). On the other hand, at high noise (amplitude $10.0$) the performance of Extragradient and Algorithm \ref{alg:cent} become comparable. In addition, we compare the performance of Alg.2 under different noise over networks with different topologies in Figure~\ref{fig:decentralized_sliding_different_noise}.

\begin{figure}[ht]
    \centering
    \includegraphics[width=1.0\textwidth]{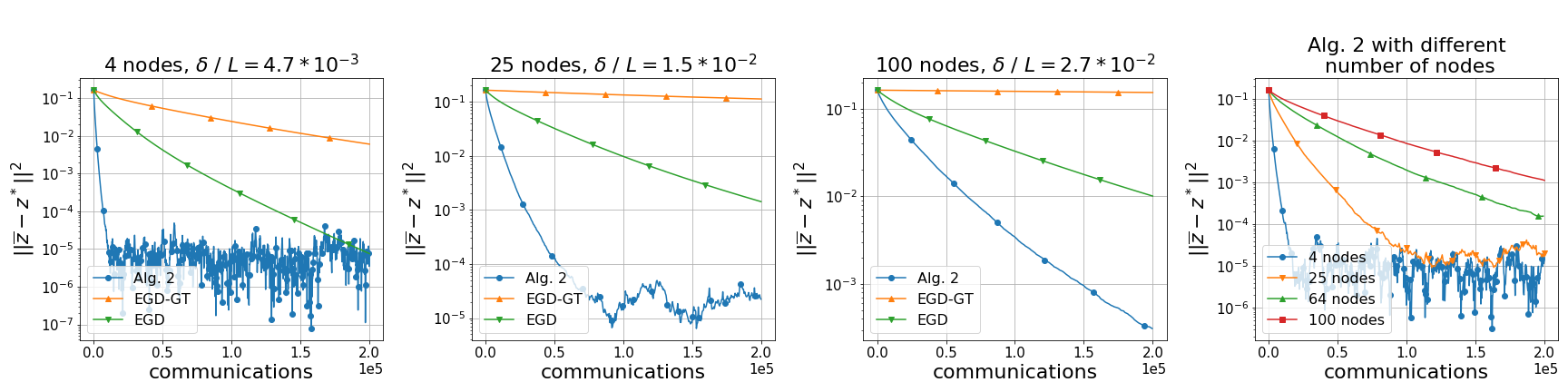}
    \caption{Decentralized case, a9a dataset, grid graph}
    \label{fig:a9a_grid}
\end{figure}
Our second experiment is using real data, specifically  LIBSVM datasets \cite{chang2011libsvm}. In this scenario, we do not use additional noise, but still can control the data similarity by choosing the number of workers. The larger  the number of workers, the less similar the local functions (less data at each node). Figure~\ref{fig:a9a_centralized} compares Algorithm \ref{alg:cent} and the Extragradient method: we plot the distance of the iterates from the solution vs. the number of communications. Quite interesting, Algorithm \ref{alg:cent} compares favorably even when the number of workers becomes large. Figure~\ref{fig:a9a_grid} compares Algorithm \ref{alg:decent} with Decentralized Extragradient method (EGD) \cite{beznosikov2021distributed} and Extragradient method with gradient-tracking (EGD-GT) \cite{9304470}. The simulations are carried out with parameters tuned according to  the    {theoretical results} in the corresponding papers.

\section{Conclusion}

We studied distristributed SPPs over networks, under data similarity. Such problems arise naturally from many applications, including machine learning and signal processing. We first derived lower complexity bounds for such problems for solution methods implementable either on star-networks or on general topologies (modeled as undirected, static graphs). These algorithms are optimal, in the sense that they achieve the lower bounds, up to log factors. The implementation of the proposed method over general network, however, is improvable: by selecting only one agent at time performing the updates, it does not fully exploit the potential computational speedup offered by the parallelism of the networking setting. Also, the use of gossip protocols to propagate the updates of a single agent across the entire network is not very efficient. Another interesting extension would be  designing  methods that take into account the asymmetry of the function $f$ with respect to the variables $x$ and $y$ (for example, various strong-convexity constants $\mu_x$ and $\mu_y$). Finally, it would be interesting to combine the proposed methods with stochastic/variance reduction techniques to alleviate the cost of local gradient computations.

\section*{Acknowledgments and Disclosure of Funding}

The research of A. Rogozin was supported by Russian Science Foundation (project No. 21-71-30005). The work of G. Scutari is supported by the  Office of Naval Research, under the Grant \# N00014-21-1-2673.
 The paper was prepared within the framework of the HSE University Basic Research Program.

\newpage
{\small
\bibliographystyle{plain}
\bibliography{lit,example_paper}
}

\newpage
\appendix

\part*{Supplementary Material}

In this appendix, we provide the proofs of the results presented in the paper; in addition to the case of   strongly-convex-strongly-concave functions (discussed therein), here we establish results also for  the case of (non strongly)  convex-concave functions.  In this latter setting,  Assumption 1 (iii) (cf.~Sec.~2) is fulfilled with $\mu = 0$; in addition, for some $G > 0$ it holds $\| F_m ( z^*) \| \leq G$, for all $m$. In the general convex-concave case, we also assume that the set $\mathcal{Z}$ is compact and introduce $\Omega$ -- the diameter of $\mathcal{Z}$.

For the sake of convenience, we summarize next the main lower/upper complexity bounds. 

\renewcommand{\arraystretch}{1.5}
\renewcommand{\tabcolsep}{10pt} 
\begin{table}[h!]
\vspace{-0.3cm}
\begin{center}
\begin{tabular}{ccc}
\hline
\multicolumn{1}{c}{} & \multicolumn{1}{c}{\textbf{lower}} & \multicolumn{1}{c}{\textbf{upper}} \\ \hline
\multicolumn{1}{c}{}   & \multicolumn{2}{c}{{\tt centralized}} \\ \hline
{\tt sc} & {$\Omega\left( \Delta\left( 1 + \frac{\delta}{\mu} \right) \log \frac{\sqn{z^{0} - z^*}}{\varepsilon} \right)$} & {$\cO \left( \Delta \left( 1 + \frac{\delta}{\mu} \right) \log \frac{\sqn{z^{0} - z^*}}{\varepsilon} \right)$} \\ \hline
{\tt c} & {$\Omega\left( \Delta \frac{\delta \Omega^2}{ \varepsilon} \right)$} & {$\cO \left( \Delta \frac{\delta \Omega^2}{ \varepsilon} \right)$} \\ \hline
\multicolumn{1}{c}{}   & \multicolumn{2}{c}{{\tt decentralized}} \\ \hline
{\tt sc} & {$\Omega\left( \frac{1}{\sqrt{\rho}}\left( 1 + \frac{\delta}{\mu} \right) \log \frac{\sqn{z^{0} - z^*}}{\varepsilon} \right)$} & {$\mathcal{ \tilde O} \left( \frac{1}{\sqrt{\rho}} \left( 1 + \frac{\delta}{\mu} \right) \log^2 \frac{\sqn{z^{0} - z^*}}{\varepsilon} \right)$} \\ \hline
{\tt c} & {$\Omega\left( \frac{1}{\sqrt{\rho}} \frac{\delta \Omega^2}{ \varepsilon} \right)$} & {$\mathcal{ \tilde O} \left( \frac{1}{\sqrt{\rho}} \frac{\delta \Omega^2}{ \varepsilon} \right)$} \\ \hline 
\end{tabular}
\vspace{0.4cm}
\caption{Comparison of lower and upper bounds on communication rounds for $\delta$-related smooth strongly-convex--strongly-concave ({\tt sc}) or convex-concave ({\tt c}) saddle-point problems in centralized and decentralized cases. Notation:
$L$ -- smothness constant of $f_m$, $\mu$ -- strongly-convex-strongly-concave constant, $\Omega$ -- diameter of optimization set, $\Delta$, $\rho$ -- diameter of communication graph and eigengap of the gossip matrix, $\varepsilon$ -- precision. In the case of upper bounds for  the convex-concave case, the convergence is in terms of the ``saddle-point residual''  [cf.~\eqref{eq:saddle-point-residual}]; for  ({\tt sc}) functions, it is  in terms of the (square) distance to the solution.}
\end{center}
\label{table1}
\end{table}

\section{Lower Complexity Bounds}

We construct the following bilinearly functions with $\delta, \mu$ and $d_x = d_y = d$. Let us consider a linear graph $G$ of $M \geq 3$ nodes. Define $p = \left\lceil\frac{M}{32}\right\rceil$; and let $B = \left\{1, \ldots p\right\}$  and  $\bar B = \{M - p + 1, \ldots, M\}$, with  $|B| = |\bar B| = p$. The distance in edges $l$ between $B$ and $\bar B$ can be bounded by $M - 2p + 1$. We then construct the following   bilinear functions on the graph:
\begin{eqnarray}
\label{t2}
f_m (x,y) = 
\begin{cases}
f_1 (x,y) = \frac{\delta}{4} x^T A_1 y + \frac{p}{M} \cdot 16\mu\|x\|^2 - \frac{p}{M} \cdot 16\mu\|y\|^2 + \frac{\delta^2}{128\mu}e_1^T y, & m \in \bar B;\\
f_2 (x,y) =\frac{\delta}{4} x^T A_2 y + \frac{p}{M} \cdot 16\mu\|x\|^2 - \frac{p}{M} \cdot 16\mu\|y\|^2, & m \in B;\\
f_3 (x,y) = \frac{p}{M} \cdot 16\mu\|x\|^2 - \frac{p}{M} \cdot 16\mu\|y\|^2, & \text{otherwise};
\end{cases}
\end{eqnarray}
where $e_1 = (1,0 \ldots, 0)$ and
\begin{eqnarray*}
A_1 = \left(
\begin{array}{cccccccc}
1&0 & & & & & &  \\
&1 &-2 & & & & &  \\
& &1 &0 & & & & \\
& & &1 &-2 & & & \\
& & & &\ldots &\ldots & & \\
& & & & &1  &-2   & \\
& & &   & & &1 &0 \\
& & &  & & & &1 \\
\end{array}
\right), ~~
A_2 = \left(
\begin{array}{cccccccc}
1&-2 & & & & & &  \\
&1 &0 & & & & &  \\
& &1 &-2 & & & & \\
& & &1 &0 & & & \\
& & & &\ldots &\ldots & & \\
& & & & &1  &0   & \\
& & &   & & &1 &-2 \\
& & &  & & & &1 \\
\end{array}
\right).
\end{eqnarray*}

Consider the global objective function:
\begin{align}
    \label{t144}
    f(x,y) &= \frac{1}{M} \sum\limits_{m=1}^M f_m(x,y) = \frac{1}{M} \left( |\bar B| \cdot f_1(x,y) + |B| \cdot f_2(x,y) + (M - |\bar B| - |B|) \cdot f_3(x,y)\right) \nonumber\\
    &= \frac{2p}{M} \cdot \frac{\delta}{4} x^T A y + \frac{p}{M} \cdot 16\mu\|x\|^2 - \frac{p}{M} \cdot 16\mu\|y\|^2 + \frac{p}{M} \cdot \frac{\delta^2}{128\mu} e_1^T y,
\end{align}
with $A = \frac{1}{2}(A_1 + A_2)$. 

It is easy to check  that
\begin{align*}
    \nabla^2_{xx} f_1 (x,y) &=\nabla^2_{xx} f_2 (x,y) = \nabla^2_{xx} f_3 (x,y) = \nabla^2_{xx} f (x,y) = \frac{p}{M} \cdot 16\mu I_x; \\
    \nabla^2_{yy} f_1 (x,y) &=\nabla^2_{yy} f_2 (x,y) = \nabla^2_{yy} f_3 (x,y) = \nabla^2_{yy} f (x,y) = \frac{p}{M} \cdot 16\mu I_y; \\
    \nabla^2_{xy} f_1 (x,y) &= \frac{\delta }{4} A_1, \quad  \nabla^2_{xy} f_2 (x,y) = \frac{\delta }{4} A_2; \\ 
    \nabla^2_{xy} f_3 (x,y) &= 0, \quad  \nabla^2_{xy} f (x,y) = \frac{2p}{M} \cdot \frac{\delta}{4} A.
\end{align*}

Note that  $f_1, f_2, f_3$ are  $L$--smooth (for $L \geq \delta$),   $\mu$-strongly-convex--strongly-concave, and  $\delta$-related; the last is a consequence of the following
$$\nabla^2_{xx} f_1 (x,y) - \nabla^2_{xx} f (x,y) = \nabla^2_{xx} f_2 (x,y) - \nabla^2_{xx} f (x,y) = \nabla^2_{xx} f_3 (x,y) - \nabla^2_{xx} f (x,y) = 0;$$
$$\nabla^2_{yy} f_1 (x,y) - \nabla^2_{yy} f (x,y) = \nabla^2_{yy} f_2 (x,y) - \nabla^2_{yy} f (x,y) = \nabla^2_{yy} f_3 (x,y) - \nabla^2_{yy} f (x,y) = 0;$$
$$\|\nabla^2_{xy} f_1 (x,y) - \nabla^2_{xy} f (x,y)\| \leq \|\nabla^2_{xy} f_1 (x,y) \| + \|  \nabla^2_{xy} f (x,y) \| \leq \delta  \left(\frac{5}{8} + \frac{p}{M} \right)  \leq \delta ;$$
$$\|\nabla^2_{xy} f_2 (x,y) - \nabla^2_{xy} f (x,y)\| \leq \|\nabla^2_{xy} f_2 (x,y) \| + \|  \nabla^2_{xy} f (x,y) \| \leq \delta  \left(\frac{5}{8} + \frac{p}{M} \right) \leq \delta ;$$
$$\|\nabla^2_{xy} f_3 (x,y) - \nabla^2_{xy} f (x,y)\| \leq \|\nabla^2_{xy} f_3 (x,y) \| + \|  \nabla^2_{xy} f (x,y) \| \leq \delta  \frac{p}{M} \leq \delta.$$

\begin{lemma} Let Problem~\eqref{t2} be solved by any method  that satisfies Definition 1. Then after $K$ communication rounds,  only the first $\left\lfloor \frac{K}{l} \right\rfloor$ coordinates of the   global output can be non-zero while  the rest of the $d-\left\lfloor \frac{K}{l} \right\rfloor$ coordinates are strictly equal to zero. Here $l = M - 2p + 1$ (distance in edges between $B$ and $\bar B$).
\end{lemma}

\textbf{Proof:} We begin introducing some notation, instrumental for our proof. Let
\begin{align*}
    E_{0} := \{ 0\}, \quad E_{K} := \text{span} \{ e_1, \ldots, e_K\}.
\end{align*}
Note that, the initialization reads   $\mathcal{M}^x_{m} = E_0$, $\mathcal{M}^y_{m} = E_0$.

Suppose that, for some $m$,  $\mathcal{M}^x_{m} = E_K$ and $\mathcal{M}^y_{m} = E_K$, at some given time. Let us analyze how $\mathcal{M}^x_{m}, \mathcal{M}^y_{m}$ can change by performing only local computations. 

Firstly, we consider the case when 
$K$ odd. We have the following: 

$\bullet$ For   machines $m$  which own $f_1$, it holds
\begin{equation*}\begin{aligned}
        \alpha x + \beta A_1 y \in \text{span} \big\{&e_1~,~ x'~,~A_1 y'~,~ A_1 A_1^T x' \big\} = E_K,\\
        \theta y - \varphi A_1^T x \in \text{span} \big\{&e_1~,~y'~,~A_1^T x' ~,~  A_1^T A_1 y' \big\} = E_K. 
\end{aligned}\end{equation*}
Since $A_1$ has a block diagonal structure with alternating blocks $1 \times 1$ and $2 \times 2$,   $A_1^{-1}$ admits the same partitions into $1 \times 1$ and $2 \times 2$ blocks on the diagonal. Therefore,   after local computations, we have  $\mathcal{M}^x_{m} = E_K$ and $\mathcal{M}^y_{m} = E_K$. The situation does not change, no matter how many local computations one  does.

$\bullet$ For machines $m$ which own $f_2$, it holds
\begin{equation*}\begin{aligned}
        \alpha x + \beta A_2 y \in \text{span} \big\{&x'~,~A_2 y'~,~ A_2 A_2^T x' \big\} = E_{K+1},\\
        \theta y - \varphi A_2^T x \in \text{span} \big\{&y'~,~A_2^T x' ~,~  A_2^T A_2 y' \big\} = E_{K+1}, 
\end{aligned}\end{equation*}
for given $x', x'' \in \mathcal{M}^x_{m}$ and  $y', y'' \in \mathcal{M}^y_{m}$. It means that, after local computations, one has  $\mathcal{M}^x_{m} = E_{K+1}$ and $\mathcal{M}^y_{m} = E_{K+1}$. Therefore,  machines with function $f_2$ can progress by one new non-zero coordinate.

This means that we constantly have to transfer progress from the group of machines with $f_1$ to the group of machines with $f_2$ and back. Initially, all devices have zero coordinates. Further, machines with $f_1$ can receive the first nonzero coordinate (but only the first, the second is not), and the rest of the devices are left with all zeros. Next, we pass the first non-zero coordinate to machines with $f_2$. To do so,    $l$ communication rounds are needed. By doing so,  they can make the second coordinate non-zero, and then transfer this progress to the machines with $f_1$. Then the process continues in the same way. This  completes the proof.\EndProof

The next lemma is   devoted to provide an approximate solution of   problem \eqref{t144}, and shows that this approximation is close to a real solution. 
The proof of the lemma follows closely that of    \cite[Lemma 3.3]{zhang2019lower}, and is reported for the sake of completeness. 

\begin{lemma}[Lemma 3.3 from \cite{zhang2019lower}]\label{lemma2}
Let $\alpha = \left(\frac{64\mu}{\delta}\right)^2$ and $q = \frac{1}{2}\left(2 + \alpha - \sqrt{\alpha^2 + 4\alpha} \right) \in (0;1)$--the smallest root of $q^2 - (2 + \alpha) q + 1 = 0$; and let define
\begin{equation*}
    \bar y^*_i = \frac{q^i}{1-q},\quad i\in [d].
\end{equation*}
The following bound holds when   $\bar y^*:=[y^*_1,\ldots y^*_d]^\top$ is used to  approximate   the solution $ y^*$: 
\begin{equation*}
    \|\bar y^* - y^*\| \leq \frac{q^{d+1}}{\alpha(1-q)}.
\end{equation*}
\end{lemma}
\textbf{Proof:} Let us write the dual function for \eqref{t144}:
\begin{equation*}
    g(y) = \frac{p}{M} \cdot \left[-\frac{1}{2}y^T \left(\frac{\delta^2}{128\mu}A^T A + 32\mu I \right)y + \frac{\delta^2}{128\mu} e_1^T y\right],
\end{equation*}
where it is not difficult to check that 
\begin{eqnarray*}
A A^T = \left(
\begin{array}{cccccccc}
1&-1 & & & & & &  \\
-1&2 &-1 & & & & &  \\
&-1 &2 & -1 & & & & \\
& & -1&2 &-1 & & & \\
& & &-1 &2 &-1 & & \\
& & & & &\ldots & & \\
& & & & &-1 &2 &-1 \\
& & & & & &-1 &2 \\
\end{array}
\right).
\end{eqnarray*}
The optimality of dual problem $\nabla g(y^*) = 0$ gives
\begin{equation*}
    \left(\frac{\delta^2}{128\mu}A^T A + 32\mu I \right)y^* = \frac{\delta^2}{128\mu} e_1,
\end{equation*}
or
\begin{equation*}
    \left(A^T A + \alpha I \right)y^* = e_1.
\end{equation*}
Equivalently, we can write 
\begin{eqnarray*}
\left\{
\begin{array}{l}
(1+\alpha)y_1^* - y_2^* = 1, \\
-y_1^* + (2 + \alpha) y^*_2 - y^*_3 = 0,\\
\ldots \\
-y_{d-2}^* + (2 + \alpha) y^*_{d-1} - y^*_d = 0,\\
-y^*_{d-1} + (2+\alpha)y^*_d = 0.
\end{array}
\right .
\end{eqnarray*}
On the other hand, the    approximation $\bar y^*$ satisfies the following set of equations:
\begin{eqnarray*}
\left\{
\begin{array}{l}
(1+\alpha)\bar y_1^* - \bar y_2^* = 1, \\
-\bar y_1^* + (2 + \alpha) \bar y^*_2 - \bar y^*_3 = 0,\\
\ldots \\
-\bar y_{d-2}^* + (2 + \alpha) \bar y^*_{d-1} - \bar y^*_d = 0,\\
-\bar y^*_{d-1} + (2+\alpha) \bar y^*_d = \frac{q^{d+1}}{1-q},
\end{array}
\right .
\end{eqnarray*}
or equivalently
\begin{equation*}
    \left(A^T A + \alpha I \right)\bar y^* = e_1 + \frac{q^{d+1}}{1-q}e_d.
\end{equation*}
Therefore,  the difference between  $\bar y^*$ and   $y^*$   reads
\begin{equation*}
    \bar y^* - y^* = \left(A^T A + \alpha I \right)^{-1}\frac{q^{d+1}}{1-q}e_d.
\end{equation*}
The statement of the lemma follow from the above equality and $\alpha^{-1} I \succeq \left(A^T A + \alpha I \right)^{-1} \succ 0$.  
\EndProof

The next lemma 
provides a lower bound for the solution of \eqref{t144} in the distributed case \eqref{t2}. The proof follows closely that of     \cite[Lemma 3.4]{zhang2019lower} and is reported for the sake of completeness. 

\begin{lemma} \label{l234}
Consider  a distributed saddle-point problem with objective function given by \eqref{t144}. For any $K$, choose any  problem size  $d \geq \max \left\{ 2 \log_q \left( \frac{\alpha}{4\sqrt{2}}\right), 2K\right\}$, where $\alpha = \left(\frac{64\mu}{\delta}\right)^2$ and $q = \frac{1}{2}\left(2 + \alpha - \sqrt{\alpha^2 + 4\alpha} \right) \in (0;1)$. Then, any output $\hat x, \hat y$ produced by  any method satisfying Definition 1 after $K$ communications rounds, is such that 
\begin{equation*}
    \|\hat x - x^*\|^2 + \|\hat y - y^*\|^2 \geq q^{\frac{2K}{l}} \frac{\| y_0 - y^*\|^2}{16}.
\end{equation*}
\end{lemma}
\textbf{Proof:} From Lemma 1 we know that after $K$ communication rounds only $k = \left\lfloor \frac{K}{l} \right\rfloor$ first coordinates in the output can be non-zero. By definition of $\bar y^*$, with $q < 1$ and $k \leq \frac{d}{2}$, we have
\begin{eqnarray*}
    \|\hat y - \bar y^*\|^2 &\geq& \sqrt{\sum\limits_{j=k+1}^n  (\bar y^*_j)^2} = \frac{q^k}{1-q} \sqrt{q^2 + q^4 + \ldots + q^{2(d-k)}} \\
    &\geq& \frac{q^k}{\sqrt{2}(1-q)} \sqrt{q^2 + q^4 + \ldots + q^{2n}} = \frac{q^k}{\sqrt{2}} \| \bar y^*\|^2 = \frac{q^k}{\sqrt{2}} \| y_0 - \bar y^*\|^2.
\end{eqnarray*}
Using Lemma \ref{lemma2} for $d \geq 2 \log_q \left(\frac{\alpha}{4\sqrt{2}} \right)$ we can guarantee that $\bar y^* \approx y^*$ (for more detailed proof see \cite{zhang2019lower}) and
\begin{equation*}
    \|\hat x - x^*\|^2 + \|\hat y - y^*\|^2 \geq \|\hat y - y^*\|^2 \geq \frac{q^{2k}}{16} \| y_0 - y^*\|^2 = q^{2\left\lfloor 
    \frac{K}{l} \right\rfloor } \frac{\| y_0 - y^*\|^2}{16} \geq q^{\frac{2K}{l}} \frac{\| y_0 - y^*\|^2}{16}.
\end{equation*}
\EndProof

\subsection{Centralized case (Theorem 1)} \label{app:low_c}
Building on the above preliminary results, we are now ready to prove our complexity lower bound as stated in Theorem 1 of the paper. The following theorem is a more detailed version of the statement in Theorem 1.

\begin{theorem}
Let $L, \mu, \delta > 0$ (with $L > \mu$ and $L > \delta$), $\Delta \in \mathcal{N}$ and $K \in \mathcal{N}$. There exists a centralized saddle-point problem on graph $\mathcal G$ for which the following statements are true:

$\bullet$ the diameter of graph $G$ is equal to $\Delta$,

$\bullet$ $f = \frac{1}{M} \sum\limits_{m=1}^M f_m : \mathcal{R}^d \times \mathcal{R}^d \to \mathcal{R}$ are $L$-Lipschitz continuous, $\mu$ -- strongly-convex-strongly-concave,

$\bullet$ $f_m$ are $L$-Lipschitz continuous, $\mu$ -- strongly-convex-strongly-concave, $\delta$-related,

$\bullet$ size $d \geq \max \left\{ 2 \log_q \left( \frac{\alpha}{4\sqrt{2}}\right), 2K\right\}$, where  $\alpha = \left(\frac{64\mu}{\delta}\right)^2$ and $q = \frac{1}{2}\left(2 + \alpha - \sqrt{\alpha^2 + 4\alpha} \right) \in (0;1)$,

$\bullet$ the solution of the problem is non-zero: $x^* \neq 0$, $y^* \neq 0$.

Then for any output $\hat z$ of any procedure (Definition 1) with $K$ communication rounds, one can obtain the following estimate:
\begin{equation*}
    \|\hat z - z^*\|^2 = \Omega\left(\exp\left( - \frac{K}{\Delta} \cdot \frac{1}{\frac{1}{8}\sqrt{1 + \left(\frac{\delta}{32\mu}\right)^2} - \frac{1}{8}} \right) \| y_0 - y^*\|^2\right).
\end{equation*}
\end{theorem}

\textbf{Proof:} It suffices to consider a linear graph with $\Delta + 1$ vertices $\{v_1, \ldots, v_{\Delta+1}\}$ and apply Lemma 1 and Lemma 3. We have
\begin{equation*}
    \left(\frac{1}{q}\right)^{\frac{2K}{l}} \geq  \frac{\| y_0 - y^*\|^2}{16(\|\hat x - x^*\|^2 + \|\hat y - y^*\|^2)}.
\end{equation*}
Taking the logarithm on both sides, we get
\begin{equation*}
    \frac{2K}{l}  \geq  \ln\left(\frac{\| y_0 - y^*\|^2}{16(\|\hat x - x^*\|^2 + \|\hat y - y^*\|^2)}\right) \frac{1}{\ln(q^{-1})}.
\end{equation*}
Next, we work with
\begin{eqnarray*}
    \frac{1}{\ln (q^{-1})} &=& \frac{1}{\ln (1+ (1-q)/q))} = \frac{1 + \frac{\alpha}{2} - \sqrt{\frac{\alpha^2}{4} + \alpha}}{\sqrt{\frac{\alpha^2}{4} + \alpha} - \frac{\alpha}{2}} \\
    &=& \frac{\sqrt{\frac{\alpha^2}{4} + \alpha} - \frac{\alpha}{2}}{\alpha} = \sqrt{\frac{1}{4} + \frac{1}{\alpha}} - \frac{1}{2} \\
    &=& \sqrt{\frac{1}{4} + \left(\frac{\delta}{64\mu}\right)^2} - \frac{1}{2}.
\end{eqnarray*}
Finally, one can then write 
\begin{equation*}
    \frac{2K}{l}  \geq  \ln\left(\frac{\| y_0 - y^*\|^2}{16(\|\hat x - x^*\|^2 + \|\hat y - y^*\|^2)}\right) \left( \frac{1}{2}\sqrt{1 + \left(\frac{\delta}{32\mu}\right)^2} - \frac{1}{2}\right),
\end{equation*}
and 
\begin{equation*}
    \exp\left(\frac{1}{\frac{1}{2}\sqrt{1 + \left(\frac{\delta}{32\mu}\right)^2} - \frac{1}{2}}\frac{2K}{l}\right)  \geq  \frac{\| y_0 - y^*\|^2}{16(\|\hat x - x^*\|^2 + \|\hat y - y^*\|^2)},
\end{equation*}
which completes the proof, with $l \geq \frac{1}{2}\Delta$.
\EndProof

\subsection{Decentralized case (Theorem 2)} \label{app:low_dec}
The lower complexity bound as stated in Theorem 2 is proved next. The next theorem is a more detailed version of Theorem 2.
\begin{theorem}
Let $L, \mu, \delta > 0$ (with $L > \mu$ and $L > \delta$), $\rho \in (0;1]$ and $K \in \mathcal{N}$. There exists a distributed saddle-point problem. For which the following statements are true:

$\bullet$ a gossip matrix $W$ have $\rho(W) = \rho$,

$\bullet$ $f = \frac{1}{M} \sum\limits_{m=1}^M f_m : \mathcal{R}^d \times \mathcal{R}^d \to \mathcal{R}$ are $L$-Lipschitz continuous, $\mu$ -- strongly-convex-strongly-concave,

$\bullet$ $f_m$ are $L$-Lipschitz continuous, $\mu$ -- strongly-convex-strongly-concave, $\delta$ - related,

$\bullet$ size $d \geq \max \left\{ 2 \log_q \left( \frac{\alpha}{4\sqrt{2}}\right), 2K\right\}$, where  $\alpha = \left(\frac{64\mu}{\delta}\right)^2$ and $q = \frac{1}{2}\left(2 + \alpha - \sqrt{\alpha^2 + 4\alpha} \right) \in (0;1)$,

$\bullet$ the solution of the problem is non-zero: $x^* \neq 0$, $y^* \neq 0$.

Then for any output $\hat z$ of any procedure (Definition 1) with $T$ communication rounds, which satisfy Definition 1, one can obtain the following estimate:
\begin{equation*}
    \|\hat z - z^*\|^2 = \Omega\left(\exp\left(\sqrt{\rho}K \cdot \frac{1}{\frac{1}{20}\sqrt{1 + \left(\frac{\delta}{32\mu}\right)^2} - \frac{1}{20}} \right)  \| y_0 - y^*\|^2\right).
\end{equation*}

\end{theorem}

\textbf{Proof:}
The proof follow similar steps as in the  proof of  \cite[Theorem 2]{scaman2017optimal}. 
Let $\gamma_M = \frac{1 - \cos \frac{\pi}{M}}{1 + \cos \frac{\pi}{M}}$ be a  decreasing sequence of positive numbers. Since $\gamma_2 = 1$ and $\lim_m \gamma_M = 0$, there exists $M \geq 2$ such that $\gamma_M \geq \rho > \gamma_{M+1}$.

$\bullet$ If $M \geq 3$, let us consider linear graph of size $M$ with vertexes $v_1, \ldots v_M$, and weighted with $w_{1,2} = 1 - a$ and  $w_{i,i+1} = 1$ for $i \geq 2$. Then we applied Lemmas 1 and 3 and get:
\begin{equation*}
    \|\hat x - x^*\|^2 + \|\hat y - y^*\|^2 \geq q^{\frac{2K}{l} } \frac{\| y_0 - y^*\|^2}{16}.
\end{equation*}
If $W_a$ is the Laplacian of the weighted graph $\mathcal G$, one can note that with $a = 0$, $\rho(W_a) = \gamma_M$, with $a = 1$, we have $\rho(W_a) = 0$. Hence, there exists $a \in (0;1]$ such that $\rho(W_a) = \rho$. Then $\rho \geq \gamma_{M+1} \geq \frac{2}{(M+1)^2}$, and $M \geq \frac{\sqrt{2}}{\sqrt{\rho}} - 1 \geq \frac{1}{4 \sqrt{\rho}}$. Finally, $l = M - 2p + 1 \geq \frac{15M}{16} - 1 \geq \frac{15}{16} \left(\frac{\sqrt{2}}{\sqrt{\rho}} - 1\right) - 1 \geq \frac{1}{5\sqrt{\rho}}$ since $\rho \leq \gamma_3 = \frac{1}{3}$. Hence,
\begin{equation*}
    \|\hat x - x^*\|^2 + \|\hat y - y^*\|^2 \geq q^{10\sqrt{\rho}K } \frac{\| y_0 - y^*\|^2}{16}.
\end{equation*}
Similarly to the proof of the previous theorem
\begin{equation}
    \label{r509}
    \exp\left(\sqrt{\rho}K \cdot \frac{1}{\frac{1}{20}\sqrt{1 + \left(\frac{\delta}{32\mu}\right)^2} - \frac{1}{20}} \right)  \geq  \frac{\| y_0 - y^*\|^2}{16(\|\hat x - x^*\|^2 + \|\hat y - y^*\|^2)}.
\end{equation}
$\bullet$ If $M = 2$, we construct a totally connected network with 3 nodes with weight $w_{1,3} = a \in [0;1]$. Let $W_a$ is the Laplacian. If $a = 0$, then the network is a linear graph and $\rho(W_a) = \gamma_3 =\frac{1}{3}$. Hence, there exists $a \in [0;1]$ such that $\rho(W_a) = \rho$. Finally, $B = \{v_1\}$, $\bar B = \{v_3\}$ and $l \geq 1 \geq \frac{1}{2\sqrt{\rho}}$. Whence it follows that in this case \eqref{r509} is also valid.
\EndProof

\subsection{Regularization and convex-concave case}
To establish  the lower bounds for the case of  (non strongly) convex-concave problems, one can use the classical trick of introducing a regularization and consider instead  the following objective function 
\begin{align*}
    g(x,y) + \frac{\varepsilon}{4\Omega^2} \cdot \|x - x^0 \|^2 - \frac{\varepsilon}{4\Omega^2}\cdot \|y - y^0 \|^2,
\end{align*}
which is   strongly-convex-strongly-concave with constant  $\mu = \frac{\varepsilon}{2\Omega^2}$, where $\varepsilon$ is a precision within the solution of the original problem is computed  and $\Omega$ is the  diameter of the  sets $\mathcal{X}$ and $\mathcal{Y}$. The resulting new SPP problem is solved to $\varepsilon/2$-precision in order to guarantee an accuracy $\varepsilon$ on the solution of  the original problem. Therefore, one can directly leverage the lower bound   estimates \eqref{cntr_b} and \eqref{decntr_b} with the   new constants above; this leads to the following  lower  bounds on the number of communications  
\begin{equation*}
    \Omega\left( \Delta \frac{\delta \Omega^2}{\varepsilon} \right),\quad \quad  \Omega\left(\frac{1}{\sqrt{\rho}} \cdot \frac{\delta \Omega^2}{\varepsilon}\right),
\end{equation*}
for the  centralized and decentralized case, respectively.

\section{Optimal algorithms}

For the general convex-concave case we introduce the following metric to measure  convergence:
\begin{align}\label{eq:saddle-point-residual}
    \text{gap}(z) = \text{gap}(x,y) := \max_{y' \in \mathcal{Y}} f(x, y') - \min_{x' \in \mathcal{X}} f(x', y).
\end{align}

\subsection{Centralized case} \label{app:op_c}

\subsubsection{Strongly-convex-strongly-concave case (Proof of   Theorem~\ref{theorem2})}
We begin  introducing some intermediate results.  Throughout this section, we tacitly subsume all the assumptions as in  Theorem~\ref{theorem2}. 
\begin{lemma} \label{lem1}
    Let $\{z^k\}$ be the sequence generated by  Algorithm \ref{alg:cent} over $\mathcal{G}$ with a master node. The following holds: 
	\begin{align}\label{eq:error-contraction}
	\sqn{z^{k+1} - z^*}
	&\leq
	(1 - \gamma \mu)\sqn{z^k - z^*} 
	- (1 - 3\gamma \mu - 4 \gamma^2 \delta^2)\sqn{z^k - \hat u^k} \nonumber
	\\&\hspace{0.4cm} 
	+ \left(2 +\frac{4 \gamma \delta^2}{\mu} + \frac{4}{\gamma \mu} + 4 \gamma^2 \delta^2\right) \sqn{u^k - \hat u^k} .
	\end{align}
\end{lemma}
\textbf{Proof:} Define $w^k = u^k + \gamma \cdot (F(z^k) - F_1(z^k) - F(u^k) + F_1(u^k))$.   Using the non-expansiveness of the Euclidean projection, we have
	\begin{align*}
		\sqn{z^{k+1} - z^*} &= \sqn{\text{proj}_{\mathcal{Z}}\left[w^{k}\right] - \text{proj}_{\mathcal{Z}}\left[z^*\right]}
		\\&\leq \sqn{w^{k} - z^*}
		\\&=
		\sqn{z^k - z^*} + 2\<w^k - z^k, z^k - z^*> + \sqn{w^k - z^k}
		\\&=
		\sqn{z^k - z^*} + 2\<w^k - z^k, \hat u^k - z^*>  + 2\<w^k - z^k,z^k - \hat u^k> + \sqn{w^k - z^k}
		\\&=
		\sqn{z^k - z^*} + 2\<w^k - z^k, \hat u^k - z^*>  + \sqn{w^k - \hat u^k} - \sqn{z^k - \hat u^k}
		\\&= 
		\sqn{z^k - z^*} + 2\< u^k + \gamma \cdot (F(z^k) - F_1(z^k) - F(u^k) + F_1(u^k)) - z^k, \hat u^k - z^*> 
		\\&\hspace{0.4cm}+
		\sqn{w^k - \hat u^k} - \sqn{z^k - \hat u^k}
		\\&= 
		\sqn{z^k - z^*} + 2\< u^k + \gamma \cdot (F(z^k) - F_1(z^k)) - z^k, \hat u^k - z^*> 
		\\&\hspace{0.4cm}-
		2\gamma \<  F(u^k) - F_1(u^k), \hat u^k - z^*> + \sqn{w^k - \hat u^k} - \sqn{z^k - \hat u^k}.
	\end{align*}
Substituting the expression of $v^k$, we have
	\begin{align*}
	\sqn{z^{k+1} - z^*}
	&\leq
	\sqn{z^k - z^*} + 2\< u^k - v^k, \hat u^k - z^*> - 2\gamma \<  F(u^k) - F_1(u^k), \hat u^k - z^*>
	\\&\hspace{0.4cm} 
	+\sqn{w^k -\hat u^k} - \sqn{z^k - \hat u^k}
	\\&= 
	\sqn{z^k - z^*} + 2\< \hat u^k - v^k, \hat u^k - z^*> - 2\gamma \<  F(u^k) - F_1(u^k), \hat u^k - z^*>
	\\&\hspace{0.4cm} 
	+ 2\< u^k - \hat u^k, \hat u^k - z^*>  +\sqn{w^k - \hat u^k} - \sqn{z^k - \hat u^k}.
	\end{align*}
Invoking  the optimality  of $\hat u^k$,   $\<\gamma F_1(\hat{u}^k) + \hat u^k - v^k, \hat u^k - z > \leq 0$ (for all $z \in \mathcal{Z}$), yields:
    \begin{align}\label{eq:delta_error_int1}
	\sqn{z^{k+1} - z^*}
	&\leq
	\sqn{z^k - z^*} - 2\gamma\< F_1(\hat{u}^k), \hat u^k - z^*> - 2\gamma \<  F(u^k) - F_1(u^k), \hat u^k - z^*>
	\nonumber\\&\hspace{0.4cm} 
	+ 2\< u^k - \hat u^k, \hat u^k - z^*> +\sqn{w^k - \hat u^k} - \sqn{z^k - \hat u^k}
	\nonumber\\&=
	\sqn{z^k - z^*} - 2\gamma\< F_1(\hat{u}^k), \hat u^k - z^*> - 2\gamma \<  F(\hat u^k) - F_1(\hat u^k), \hat u^k - z^*>
	\nonumber\\&\hspace{0.4cm} 
	+ 2\< \gamma (F(\hat u^k) - F_1(\hat u^k) - F(u^k) + F_1(u^k)) + u^k - \hat u^k, \hat u^k - z^*> \nonumber
	\\&\hspace{0.4cm} 
	+\sqn{w^k - \hat u^k} - \sqn{z^k - \hat u^k}.
	\end{align}
Invoking the optimality  of the solution $z^*$: $\<\gamma F(z^*), z^* - z > \leq 0$ (for all $z \in \mathcal{Z}$) along with the $\mu$-strong convexity-strong concavity of $f$, we obtain
    \begin{align*}
	\sqn{z^{k+1} - z^*}
	&\leq
	\sqn{z^k - z^*} - 2\gamma\< F(\hat{u}^k) - F(z^*), \hat u^k - z^*> 
	\\&\hspace{0.4cm} 
	+ 2\< \gamma (F(\hat u^k) - F_1(\hat u^k) - F(u^k) + F_1(u^k)) + u^k - \hat u^k, \hat u^k - z^*> 
	\\&\hspace{0.4cm} 
	+\sqn{w^k - \hat u^k} - \sqn{z^k - \hat u^k}
	\\&\leq
	\sqn{z^k - z^*} - 2\gamma \mu \sqn{\hat u^k - z^*} 
	\\&\hspace{0.4cm} 
	+ 2\< \gamma (F(\hat u^k) - F_1(\hat u^k) - F(u^k) + F_1(u^k)) + u^k - \hat u^k, \hat u^k - z^*> 
	\\&\hspace{0.4cm} 
	+\sqn{w^k - \hat u^k} - \sqn{z^k - \hat u^k}.
	\end{align*}
By Young's inequality, we have
	\begin{align*}
	\sqn{z^{k+1} - z^*}
	&\leq
	\sqn{z^k - z^*} - 2\gamma \mu \sqn{\hat u^k - z^*} 
	\\&\hspace{0.4cm} 
	+ \frac{2 }{\gamma \mu} \sqn{\gamma (F(\hat u^k) - F_1(\hat u^k) - F(u^k) + F_1(u^k)) + u^k - \hat u^k} + \frac{\gamma \mu}{2}\sqn{ \hat u^k - z^*} 
	\\&\hspace{0.4cm} 
	+\sqn{w^k - \hat u^k} - \sqn{z^k - \hat u^k}
	\\&\leq 
	\sqn{z^k - z^*} - \frac{3\gamma \mu}{2}\sqn{\hat u^k - z^*} 
	\\&\hspace{0.4cm} 
	+ \frac{4 \gamma}{\mu} \sqn{F(\hat u^k) - F_1(\hat u^k) - F(u^k) + F_1(u^k)} +  \frac{4}{\gamma \mu} \sqn{u^k - \hat u^k}  
	\\&\hspace{0.4cm} 
	+\sqn{u^k + \gamma \cdot (F(z^k) - F_1(z^k) - F(u^k) + F_1(u^k)) - \hat u^k} - \sqn{z^k - \hat u^k}
	\\&= 
	\sqn{z^k - z^*} - \frac{3\gamma \mu}{2}\sqn{\hat u^k - z^*} 
	\\&\hspace{0.4cm} 
	+ \frac{4 \gamma}{\mu} \sqn{F(\hat u^k) - F_1(\hat u^k) - F(u^k) + F_1(u^k)} +  \frac{4}{\gamma \mu} \sqn{u^k - \hat u^k}  
	\\&\hspace{0.4cm} 
	+2\sqn{u^k - \hat u^k} +2 \gamma^2\sqn{F(z^k) - F_1(z^k) - F(u^k) + F_1(u^k)} - \sqn{z^k - \hat u^k}.
	\end{align*}
Note that the function $f-f_1$ is $\delta$-smooth, since $\sqn{\nabla_{xx} f- \nabla_{xx}f_1} \leq \delta$, $\sqn{\nabla_{xy} f- \nabla_{xy}f_1} \leq \delta$, $\sqn{\nabla_{yy} f- \nabla_{yy}f_1} \leq \delta$; therefore,
    \begin{align*}
	\sqn{z^{k+1} - z^*}
	&\leq
	\sqn{z^k - z^*} - \frac{3\gamma \mu}{2}\sqn{\hat u^k - z^*} 
	\\&\hspace{0.4cm} 
	+ \frac{4 \gamma \delta^2}{\mu} \sqn{u^k - \hat u^k} +  \frac{4}{\gamma \mu} \sqn{u^k - \hat u^k}  
	\\&\hspace{0.4cm} 
	+2\sqn{u^k - \hat u^k} +2 \gamma^2 \delta^2 \sqn{z^k - u^k} - \sqn{z^k - \hat u^k}
	\\&\leq 
	\sqn{z^k - z^*} - \frac{3\gamma \mu}{2}\sqn{\hat u^k - z^*} - (1 - 4 \gamma^2 \delta^2)\sqn{z^k - \hat u^k}
	\\&\hspace{0.4cm} 
	+ \left(2 +\frac{4 \gamma \delta^2}{\mu} + \frac{4 }{\gamma \mu} + 4 \gamma^2 \delta^2\right) \sqn{u^k - \hat u^k} .
	\end{align*}
	Finally, using   $\sqn{a+b} \geq \frac{2}{3}\sqn{a} - 2\sqn{b}$, we obtain the desired result \eqref{eq:error-contraction}.
\EndProof

\begin{theorem}\label{th:inner-outer-convergence}
Let $\{z^k\}$ the sequence generated by Algorithm \ref{alg:cent} (in the setting of  Theorem~\ref{theorem2}) with the step-size $\gamma$ given by \begin{align}
\label{gamma1}
    \gamma = \min\left\{ \frac{1}{12 \mu}, \frac{1}{4 \delta}\right\}.
\end{align}
Let each subproblem \eqref{inter} be solved up to (relative) precision $\tilde e$,   
\begin{align}
\label{e1}
    \tilde e = \frac{1}{2\left(2 +\frac{4 \gamma \delta^2}{\mu} + \frac{4}{\gamma \mu} + 4 \gamma^2 \delta^2\right)}.
\end{align}
Then,   $\sqn{ z^{K} - z^*} \leq \varepsilon$ after 
\begin{align}
\label{K1}
    K = \cO \left( \left( 1 + \frac{\delta}{\mu} \right) \log \frac{\sqn{z^{0} - z^*}}{\varepsilon} \right)
    \quad 
     \text{iterations/communications}.
\end{align}
\end{theorem}

\textbf{Proof:} The output $u^k$ produced by inner method satisfies 
	\begin{equation*}
		\sqn{u^k - \hat{u}^k}\leq \tilde e \sqn{z^k - \hat u^k}.
	\end{equation*}
Combining this fact  and Lemma \ref{lem1} yields 
\begin{align*}
	\sqn{z^{k+1} - z^*} &\leq
	(1 - \gamma \mu)\sqn{z^k - z^*} 
	- (1 - 3\gamma \mu - 4 \gamma^2 \delta^2)\sqn{z^k - \hat u^k}
	\\&\hspace{0.4cm} 
	+ \left(2 +\frac{4 \gamma \delta^2}{\mu} + \frac{4}{\gamma \mu} + 4 \gamma^2 \delta^2\right) \tilde e \sqn{z^k - \hat u^k}\\
	& \overset{\eqref{e1}}{\leq} 
	(1 - \gamma \mu)\sqn{z^k - z^*} 
	- \left(\frac{1}{2} - 3\gamma \mu - 4 \gamma^2 \delta^2\right)\sqn{z^k - \hat u^k}.
\end{align*}
The proof is completed by choosing $\gamma$ according to  \eqref{gamma1}.
\EndProof

\begin{corollary}
\label{lem2}
Let we solve the subproblem   \eqref{inter} via Extragradient method with starting point $z^k$ and
	\begin{equation}
	\label{T1}
		T = \cO  \left(\left(1 + \gamma L\right)\log\frac{1}{\tilde e} \right)
	\end{equation}
iterations.	Then we can estimate the total number local iterations at the server side by 
\begin{equation*}
		\cO \left(  \left(1 + \frac{\delta}{\mu} + \frac{L}{\mu}\right)\log\frac{1}{\tilde e} \log \frac{\sqn{z^{0} - z^*}}{\varepsilon} \right).
	\end{equation*}
\end{corollary} 
\textbf{Proof:} Firstly, one can note that after $T$ iterations of Extragradient method from \eqref{T1} we can achieve $\tilde e$ precision. It follows readily from the convergence of  Extragradient method \cite{beznosikov2021distributed} and the fact that the objective function  in \eqref{inter} is $1$-strongly-convex-strongly-concave and  $(1 + \gamma L)$-smooth. 
Then we can estimate the total number of local iterations at the  server side, namely:
\begin{align*}
    K \cdot T  &= 
    \cO \left( \frac{1}{\gamma \mu} \left(1 + \gamma L\right)\log\frac{1}{\tilde e} \log \frac{\sqn{z^{0} - z^*}}{\varepsilon} \right)
    \\&= \cO \left(  \left(\frac{1}{\gamma \mu} + \frac{L}{\mu}\right)\log\frac{1}{\tilde e} \log \frac{\sqn{z^{0} - z^*}}{\varepsilon} \right)
    \\&= \cO \left(  \left(1 + \frac{\delta}{\mu} + \frac{L}{\mu}\right)\log\frac{1}{\tilde e} \log \frac{\sqn{z^{0} - z^*}}{\varepsilon} \right).
\end{align*}
\EndProof

  
\textbf{Remark.} 
If the server is located in the center of a graph with a diameter $\Delta$, then an additional factor $\Delta$ will appear in the total number of communications \eqref{K1}. 


\subsubsection{Convex-Concave case}

\begin{lemma}
    For one iteration of Algorithm \ref{alg:cent}, the following estimate holds:
	\begin{align}
	\label{t808}
	2\gamma\< F({u}^k),  u^k - z>
    &\leq
	\sqn{z^k - z} - \sqn{z^{k+1} - z}  - ( 1 - 2 \gamma^2 \delta^2 ) \sqn{z^k - u^k} \notag
	\\&\hspace{0.4cm} 
	+ \left( 8\gamma L \Omega + 6\gamma G + 2 \Omega\right) \|{u}^k - \hat{u}^k\| +2\sqn{u^k - \hat u^k}  .
	\end{align}
\end{lemma}
\textbf{Proof:} The proof follows similar steps as that of Lemma \ref{lem1}, with the difference that $z^*$ therein is replaced here with any $z \in \mathcal{Z}$. Specifically,  recalling the first equality in  \eqref{eq:delta_error_int1}, we have 
    \begin{align*}
	\sqn{z^{k+1} - z}
	&\leq
	\sqn{z^k - z} - 2\gamma\< F_1(\hat{u}^k), \hat u^k - z> - 2\gamma \<  F(u^k) - F_1(u^k), \hat u^k - z>
	\\&\hspace{0.4cm} 
	+ 2\< u^k - \hat u^k, \hat u^k - z> +\sqn{w^k - \hat u^k} - \sqn{z^k - \hat u^k}
	\\&= \sqn{z^k - z} - 2\gamma\< F_1({u}^k), u^k - z> - 2\gamma \<  F(u^k) - F_1(u^k), u^k - z>
	\\&\hspace{0.4cm} 
	+ 2\gamma\< F_1({u}^k) - F_1(\hat{u}^k), u^k - z> +2\gamma\< F_1(\hat{u}^k),  u^k - \hat u^k> 
	\\&\hspace{0.4cm} 
	+2\gamma \<  F(u^k) - F_1(u^k),  u^k - \hat u^k> + 2\< u^k - \hat u^k, \hat u^k - z> 
	\\&\hspace{0.4cm} 
	+\sqn{w^k - \hat u^k} - \sqn{z^k - \hat u^k}.
	\end{align*}
Small rearrangement gives
    \begin{align*}
	2\gamma\< F({u}^k),  u^k - z>
	&\leq
	\sqn{z^k - z} - \sqn{z^{k+1} - z} 
	\\&\hspace{0.4cm} 
	+ 2\gamma\< F_1({u}^k) - F_1(\hat{u}^k), u^k - z> +2\gamma\< F_1(\hat{u}^k),  u^k - \hat u^k> 
	\\&\hspace{0.4cm} 
	+2\gamma \<  F(u^k) - F_1(u^k),  u^k - \hat u^k> + 2\< u^k - \hat u^k, \hat u^k - z> 
	\\&\hspace{0.4cm} 
	+\sqn{w^k - \hat u^k} - \sqn{z^k - \hat u^k}
	\\&\leq
	\sqn{z^k - z} - \sqn{z^{k+1} - z} 
	\\&\hspace{0.4cm} 
	+ 2\gamma \| F_1({u}^k) - F_1(\hat{u}^k)\| \cdot  \| u^k - z\| +2\gamma\| F_1(\hat{u}^k)\| \cdot \|  u^k - \hat u^k\| 
	\\&\hspace{0.4cm} 
	+2\gamma \| F(u^k) - F_1(u^k) \| \cdot  \| u^k - \hat u^k \| + 2 \| u^k - \hat u^k \| \cdot \| \hat u^k - z \| 
	\\&\hspace{0.4cm} 
	+\sqn{w^k - \hat u^k} - \sqn{z^k - \hat u^k}.
	\end{align*}
Invoking the definition of $w^k = u^k + \gamma \cdot (F(z^k) - F_1(z^k) - F(u^k) + F_1(u^k))$, we get
    \begin{align*}
	2\gamma\< F({u}^k),  u^k - z>
    &\leq
	\sqn{z^k - z} - \sqn{z^{k+1} - z} 
	\\&\hspace{0.4cm} 
	+ 2\gamma \| F_1({u}^k) - F_1(\hat{u}^k)\| \cdot  \| u^k - z\| +2\gamma\| F_1(\hat{u}^k)\| \cdot \|  u^k - \hat u^k\| 
	\\&\hspace{0.4cm} 
	+2\gamma \| F(u^k) - F_1(u^k) \| \cdot  \| u^k - \hat u^k \| + 2 \| u^k - \hat u^k \| \cdot \| \hat u^k - z \| 
	\\&\hspace{0.4cm} 
	+\sqn{u^k + \gamma \cdot (F(z^k) - F_1(z^k) - F(u^k) + F_1(u^k)) - \hat u^k} - \sqn{z^k - \hat u^k}
	\\&\leq
	\sqn{z^k - z} - \sqn{z^{k+1} - z} 
	\\&\hspace{0.4cm} 
	+ 2\gamma \| F_1({u}^k) - F_1(\hat{u}^k)\| \cdot  \| u^k - z\| +2\gamma\| F_1(\hat{u}^k)\| \cdot \|  u^k - \hat u^k\| 
	\\&\hspace{0.4cm} 
	+2\gamma \| F(u^k) - F_1(u^k) \| \cdot  \| u^k - \hat u^k \| + 2 \| u^k - \hat u^k \| \cdot \| \hat u^k - z \| 
	\\&\hspace{0.4cm} 
	+2\sqn{u^k - \hat u^k} + 2 \gamma^2 \sqn{  F(z^k) - F_1(z^k) - F(u^k) + F_1(u^k)} 
	\\&\hspace{0.4cm} - \sqn{z^k - \hat u^k}.
	\end{align*}
Then we use smoothness of $f -f_1$, $f$, $f_1$ and obtain 
    \begin{align*}
	2\gamma\< F({u}^k),  u^k - z>
    &\leq
	\sqn{z^k - z} - \sqn{z^{k+1} - z} 
	\\&\hspace{0.4cm} 
	+ 2\gamma L \|{u}^k - \hat{u}^k\| \cdot  \Omega +2\gamma(G + L \Omega) \cdot \|  u^k - \hat u^k\| 
	\\&\hspace{0.4cm} 
	+4\gamma (G + L \Omega) \| \cdot  \| u^k - \hat u^k \| + 2 \Omega \cdot \| u^k - \hat u^k \|  
	\\&\hspace{0.4cm} 
	+2\sqn{u^k - \hat u^k} + 2 \gamma^2 \delta^2 \sqn{z^k - u^k} - \sqn{z^k - \hat u^k}
	\\&=
	\sqn{z^k - z} - \sqn{z^{k+1} - z}  - ( 1 - 2 \gamma^2 \delta^2 ) \sqn{z^k - u^k}
	\\&\hspace{0.4cm} 
	+ \left( 8\gamma L \Omega + 6\gamma G + 2 \Omega\right) \|{u}^k - \hat{u}^k\| +2\sqn{u^k - \hat u^k}  .
	\end{align*}
\EndProof
Here we additionally used the diameter $\Omega$ of $\mathcal{Z}$ and simple fact: 
    \begin{align}
    \label{t101}
    \| F_1(\hat{u}^k)\| - G \leq \| F_1(\hat{u}^k)\| - \|F_1(z^*) \| \leq \| F_1(\hat{u}^k)- F_1(z^*) \| \leq L \Omega.
    \end{align}

\begin{theorem}
Let problem \eqref{inter} be solved by Extragradient with precision $e$:
\begin{align}
\label{e2}
    e = \min \left\{ \frac{\varepsilon}{\delta} ; \frac{\varepsilon^2}{(L \Omega + G + \delta \Omega)^2} \right\}
\end{align}
and number of iterations $T$:
\begin{align*}
    T =  \cO  \left(\left(1 + \gamma L\right)\log\frac{\Omega^2}{e} \right).
\end{align*}
Additionally, let us choose stepsize $\gamma$ as follows
\begin{align}
\label{gamma2}
    \gamma = \frac{1}{2\delta}.
\end{align}
Then it holds that $\text{gap}(z^{K}_{avg}) \sim \varepsilon$ after 
\begin{align}
\label{K2}
    K = \cO \left( \frac{\delta \Omega^2}{ \varepsilon} \right)\quad \text{iterations},
\end{align}
where $z^{K}_{avg}$ define as follows: $x^{K}_{avg} = \frac{1}{K}\sum_{k=0}^K u^{k}_x$, $y^{K}_{avg} = \frac{1}{K}\sum_{k=0}^K u^{k}_y$.
\end{theorem}

\textbf{Proof:} Summing \eqref{t808} over all $k$ from $0$ to $K$
\begin{align*}
	2\gamma \sum\limits_{k=0}^K\< F({u}^k),  u^k - z> &\leq \sqn{z^0 - z} - ( 1 - 2 \gamma^2 \delta^2 ) \sum\limits_{k=0}^K \sqn{z^k - u^k}
	\\&\hspace{0.4cm} 
	+ \left( 8\gamma L \Omega + 6\gamma G + 2 \Omega\right) \sum\limits_{k=0}^K \|{u}^k - \hat{u}^k\| +2 \sum\limits_{k=0}^K \sqn{u^k - \hat u^k} .
\end{align*}
Then, by $x^{K}_{avg} = \frac{1}{K}\sum_{k=0}^K u^{k}_x$ and $y^{K}_{avg} = \frac{1}{K}\sum_{k=0}^K u^{k}_y$, Jensen's inequality and convexity-concavity of $f$:
\begin{align*}
    \text{gap}(z^{K}_{avg})
    &\leq \max\limits_{y' \in \mathcal{Y}} f\left(\frac{1}{K} \left(\sum^K_{k=0} u^{k}_x \right), y'\right) - \min\limits_{x' \in \mathcal{X}} f\left(x', \frac{1}{K} \left(\sum^K_{k=0} u^{k}_y \right)\right) 
    \nonumber \\
    &\leq \max\limits_{y' \in \mathcal{Y}} \frac{1}{K} \sum^K_{k=0} f(u^{k}_x, y')  - \min\limits_{x' \in \mathcal{X}} \frac{1}{K} \sum^K_{k=0} f(x', u^{k}_y).
\end{align*}
Given the fact of linear independence of $x'$ and $y'$:
\begin{align*}
    \text{gap}(z^{K}_{avg}) &\leq \max\limits_{(x', y') \in \mathcal{Z}}\frac{1}{K} \sum^K_{k=0}  \left(f(x^{K}, y')  - f(x', u^{k}_y) \right).
\end{align*}
Using convexity and concavity of the function $f$:
\begin{align*}
    \text{gap}(z^{K}_{avg}) &\leq  \max\limits_{(x', y') \in \mathcal{Z}}\frac{1}{K} \sum^K_{k=0}  \left(f(u^{k}_x, y')  - f(x', u^{k}_y) \right)  \nonumber \\
    &= \max\limits_{(x', y') \in \mathcal{Z}} \frac{1}{K} \sum^K_{k=0} \left(f(u^{k}_x, y') -f(u^{k}_x, u^{k}_y) + f(u^{k}_x, u^{k}_y) - f(x', u^{k}_y) \right) \nonumber \\
    &\leq \max\limits_{(x', y') \in \mathcal{Z}} \frac{1}{K} \sum^K_{k=0} \left(\langle \nabla_y f (u^{k}_x, u^{k}_y), y'-u^{k}_y \rangle + \langle \nabla_x f (u^{k}_x, u^{k}_y), u^{k}_x-x' \rangle \right) \nonumber \\
    &\leq \max\limits_{z \in \mathcal{Z}} \frac{1}{K} \sum^K_{k=0} \langle  F(u^{k}), u^{k} - z\rangle.
\end{align*}
Then it gives with our choice of $\gamma$
\begin{align*}
	\text{gap}(z^{K}_{avg}) &\leq \max\limits_{z \in \mathcal{Z}} \frac{\sqn{z^0 - z}}{2\gamma K} + \frac{\left( 4\gamma L \Omega + 3\gamma G + \Omega\right)}{\gamma K} \sum\limits_{k=0}^K \|{u}^k - \hat{u}^k\| + \frac{1}{\gamma K} \sum\limits_{k=0}^K \sqn{u^k - \hat u^k} 
	\\&\leq \frac{\Omega^2}{2\gamma K} + \left( 4L \Omega + 3G + \frac{\Omega}{\gamma}\right) \sqrt{e} + \frac{1}{\gamma} e 
	\\&= \frac{\delta \Omega^2}{ K} + \left( 4L \Omega + 3G + 2 \delta \Omega\right) \sqrt{e} + 2\delta e .
\end{align*}
$e$ from \eqref{e2} is completed the proof.
\EndProof

\textbf{Remark.} \eqref{K2} also corresponds to the number of communication rounds. It is also easy to estimate the total number of local iterations on server:
\begin{align*}
    K \times T  &= 
    \cO \left( \frac{\delta \Omega^2}{ \varepsilon} \left(1 + \gamma L\right)\log\frac{\Omega^2}{e} \right)
    \\&= \cO \left(  \frac{\delta \Omega^2}{ \varepsilon} \left(1 + \frac{L}{\delta}\right)\log\frac{\Omega^2}{e}\right)
    \\&= \cO \left(  \frac{(L + \delta) \Omega^2}{ \varepsilon}\log\frac{\Omega^2}{e}\right).
\end{align*}

	
\subsection{Decentralized case} \label{app:op_dec}

Before moving on to the proofs of the decentralized case, let us understand the {\tt AccGossip} convergence  \cite{liu2011accelerated,ye2020multi}:

\begin{lemma} \label{l8} Assume that $\{y_m\}^M_{m=1}$ are output of Algorithm \ref{accgossip} with input $\{ x_m\}^M_{m=1}$. Then it holds that
\begin{align}
\label{400}
\sum\limits_{m=1}^M\|y_m - \bar y\|^2 \leq \left(1 - \sqrt{\rho}\right)^{2H} \left(\sum\limits_{m=1}^M\|x_m - \bar x\|^2\right).
\end{align}
And $\bar x = \frac{1}{M} \sum_{m=1}^M x_m = \frac{1}{M} \sum_{m=1}^M y_m =  \bar y$.
\end{lemma}
From this lemma it holds that for any $i$
\begin{align}
\label{401}
    \|y_i - \bar y\|^2 \leq \left(1 - \sqrt{\rho}\right)^{2H} \left(\sum\limits_{m=1}^M\|x_m - \bar x\|^2\right).
\end{align}
and 
\begin{align}
\label{402}
\|y_i - \bar y\| \leq \left(1 - \sqrt{\rho}\right)^{H} \sqrt{\left(\sum\limits_{m=1}^M\|x_m - \bar x\|^2\right)}.
\end{align}

\subsubsection{Strongly-convex-strongly-concave case}

\begin{lemma} \label{lem4}
    For one iteration of Algorithm 2, the following estimate holds:
	\begin{align*}
	\sqn{z^{k+1}_{m_k}  - z^*}
	&\leq (1 - \gamma \mu)\sqn{z^{k}_{m_k} - z^*} - (1 - 3\gamma \mu - 12\gamma^2\delta^2)\sqn{z^{k}_{m_k} - \hat u^{k}_{m_k}}
    \\&\hspace{0.4cm} + \left(2 + 12\gamma^2\delta^2+ \frac{4}{\gamma \mu } + \frac{8\gamma \delta^2}{\mu}\right)\sqn{ \tilde u^k_{m_k}- \hat u^k_{m_k}} 
    \\&\hspace{0.4cm} + 6\gamma^2\sqn{ \bar F^k_{m_k} - F(z^k_{m_k})} + \left(6\gamma^2 + \frac{8\gamma}{\mu}\right)\sqn{ \bar F^{k+1/2}_{m_k} - F(\tilde u^k_{m_k})}
	\\&\hspace{0.4cm} + 2\< \hat z^{k+1}_{m_k} - \tilde z^{k+1}_{m_k}, \tilde z^{k+1}_{m_k} - z^* > + \sqn{\tilde z^{k+1}_{m_k} - \hat z^{k+1}_{m_k}}.
	\end{align*}
\end{lemma}
\textbf{Proof:} Using non-expansiveness of the Euclidean projection, we get
	\begin{align*}
		\sqn{z^{k+1}_{m_k} - z^*} &= \sqn{\text{proj}_{\mathcal{Z}}\left[\hat z^{k+1}_{m_k}\right] - \text{proj}_{\mathcal{Z}}\left[z^*\right]}
		\\&\leq \sqn{\hat z^{k+1}_{m_k} - z^*}
	    \\&= \sqn{\tilde z^{k+1}_{m_k} - z^*} + 2\< \hat z^{k+1}_{m_k} - \tilde z^{k+1}_{m_k}, \tilde z^{k+1}_{m_k} - z^* > + \sqn{\tilde z^{k+1}_{m_k} - \hat z^{k+1}_{m_k}}
		\\&=
		\sqn{z^{k}_{m_k} - z^*} + 2\<\tilde z^{k+1}_{m_k} - z^{k}_{m_k}, z^{k}_{m_k} - z^*> + \sqn{\tilde z^{k+1}_{m_k} - z^{k}_{m_k}}
		\\&\hspace{0.4cm} + 2\< \hat z^{k+1}_{m_k} - \tilde z^{k+1}_{m_k}, \tilde z^{k+1}_{m_k} - z^* > + \sqn{\tilde z^{k+1}_{m_k} - \hat z^{k+1}_{m_k}}
		\\&=
		\sqn{z^{k}_{m_k} - z^*} + 2\<\tilde z^{k+1}_{m_k} - z^{k}_{m_k}, \hat u^{k}_{m_k} - z^*> + 2\<\tilde z^{k+1}_{m_k} - z^{k}_{m_k}, z^{k}_{m_k} - \hat u^{k}_{m_k}> 
		\\&\hspace{0.4cm} + \sqn{\tilde z^{k+1}_{m_k} - z^{k}_{m_k}} + 2\< \hat z^{k+1}_{m_k} - \tilde z^{k+1}_{m_k}, \tilde z^{k+1}_{m_k} - z^* > + \sqn{\tilde z^{k+1}_{m_k} - \hat z^{k+1}_{m_k}}
		\\&=
		\sqn{z^{k}_{m_k} - z^*} + 2\<\tilde z^{k+1}_{m_k} - z^{k}_{m_k}, \hat u^{k}_{m_k} - z^*> + \sqn{\tilde z^{k+1}_{m_k} - \hat u^{k}_{m_k}} - \sqn{z^{k}_{m_k} - \hat u^{k}_{m_k}}
		\\&\hspace{0.4cm} + 2\< \hat z^{k+1}_{m_k} - \tilde z^{k+1}_{m_k}, \tilde z^{k+1}_{m_k} - z^* > + \sqn{\tilde z^{k+1}_{m_k} - \hat z^{k+1}_{m_k}}
		\\&= 
		\sqn{z^{k}_{m_k} - z^*} 
		\\&\hspace{0.4cm}
		+ 2\<\tilde u^k_{m_k} + \gamma \cdot (\bar F^k_{m_k} - F_{m_k} (z^k_{m_k}) - \bar F^{k+1/2}_{m_k} + F_{m_k} (\tilde u^k_{m_k})) - z^{k}_{m_k}, \hat u^{k}_{m_k} - z^*> 
	    \\&\hspace{0.4cm} + \sqn{\tilde z^{k+1}_{m_k} - \hat u^{k}_{m_k}} - \sqn{z^{k}_{m_k} - \hat u^{k}_{m_k}} 
	    \\&\hspace{0.4cm}
	    + 2\< \hat z^{k+1}_{m_k} - \tilde z^{k+1}_{m_k}, \tilde z^{k+1}_{m_k} - z^* > + \sqn{\tilde z^{k+1}_{m_k} - \hat z^{k+1}_{m_k}}
		\\&= 
		\sqn{z^{k}_{m_k} - z^*} 
		\\&\hspace{0.4cm} + 2\<\tilde u^k_{m_k} + \gamma \cdot (\bar F^k_{m_k} - F_{m_k} (z^k_{m_k})) - z^{k}_{m_k}, \hat u^{k}_{m_k} - z^*> 
	    \\&\hspace{0.4cm} - 2\gamma\< \bar F^{k+1/2}_{m_k} - F_{m_k} (\tilde u^k_{m_k})), \hat u^{k}_{m_k} - z^*> + \sqn{\tilde z^{k+1}_{m_k} - \hat u^{k}_{m_k}} - \sqn{z^{k}_{m_k} - \hat u^{k}_{m_k}}
		\\&\hspace{0.4cm} + 2\< \hat z^{k+1}_{m_k} - \tilde z^{k+1}_{m_k}, \tilde z^{k+1}_{m_k} - z^* > + \sqn{\tilde z^{k+1}_{m_k} - \hat z^{k+1}_{m_k}}.
	\end{align*}
Substituting the expression for $v^k_{m_k}$, we have
	\begin{align*}
	\sqn{z^{k+1}_{m_k}  - z^*}
	&\leq
	\sqn{z^{k}_{m_k} - z^*} + 2\<\tilde u^k_{m_k} - v^k_{m_k}, \hat u^{k}_{m_k} - z^*> 
    \\&\hspace{0.4cm} - 2\gamma\< \bar F^{k+1/2}_{m_k} - F_{m_k} (\tilde u^k_{m_k})), \hat u^{k}_{m_k} - z^*> + \sqn{\tilde z^{k+1}_{m_k} - \hat u^{k}_{m_k}} - \sqn{z^{k}_{m_k} - \hat u^{k}_{m_k}}
	\\&\hspace{0.4cm} + 2\< \hat z^{k+1}_{m_k} - \tilde z^{k+1}_{m_k}, \tilde z^{k+1}_{m_k} - z^* > + \sqn{\tilde z^{k+1}_{m_k} - \hat z^{k+1}_{m_k}}
	\\&= 
	\sqn{z^{k}_{m_k} - z^*} + 2\<\hat u^{k}_{m_k} - v^k_{m_k}, \hat u^{k}_{m_k} - z^*> 
    \\&\hspace{0.4cm} - 2\gamma\< \bar F^{k+1/2}_{m_k} - F_{m_k} (\tilde u^k_{m_k})), \hat u^{k}_{m_k} - z^*> + \sqn{\tilde z^{k+1}_{m_k} - \hat u^{k}_{m_k}} - \sqn{z^{k}_{m_k} - \hat u^{k}_{m_k}}
	\\&\hspace{0.4cm} + 2\<\tilde u^k_{m_k} - \hat u^{k}_{m_k}, \hat u^{k}_{m_k} - z^*>  + 2\< \hat z^{k+1}_{m_k} - \tilde z^{k+1}_{m_k}, \tilde z^{k+1}_{m_k} - z^* > + \sqn{\tilde z^{k+1}_{m_k} - \hat z^{k+1}_{m_k}}.
	\end{align*}
According to the optimal condition for $\hat u^{k}_{m_k}$: $\<\gamma F_{m_k}(\hat u^{k}_{m_k}) + \hat u^{k}_{m_k} - v^k_{m_k}, \hat u^{k}_{m_k} - z > \leq 0$ (for all $z \in \mathcal{Z}$),
    \begin{align*}
	\sqn{z^{k+1}_{m_k}  - z^*}
	&\leq
	\sqn{z^{k}_{m_k} - z^*} - 2\gamma\<F_{m_k}(\hat u^{k}_{m_k}), \hat u^{k}_{m_k} - z^*> 
    \\&\hspace{0.4cm} - 2\gamma\< \bar F^{k+1/2}_{m_k} - F_{m_k} (\tilde u^k_{m_k}), \hat u^{k}_{m_k} - z^*> + \sqn{\tilde z^{k+1}_{m_k} - \hat u^{k}_{m_k}} - \sqn{z^{k}_{m_k} - \hat u^{k}_{m_k}}
	\\&\hspace{0.4cm} + 2\<\tilde u^k_{m_k} - \hat u^{k}_{m_k}, \hat u^{k}_{m_k} - z^*>  + 2\< \hat z^{k+1}_{m_k} - \tilde z^{k+1}_{m_k}, \tilde z^{k+1}_{m_k} - z^* > + \sqn{\tilde z^{k+1}_{m_k} - \hat z^{k+1}_{m_k}}
	\\&=
	\sqn{z^{k}_{m_k} - z^*} 
	- 2\gamma\<F_{m_k}(\hat u^{k}_{m_k}), \hat u^{k}_{m_k} - z^*> 
	\\&\hspace{0.4cm}- 2\gamma\< F(\hat u^k_{m_k}) - F_{m_k} (\hat u^k_{m_k}), \hat u^{k}_{m_k} - z^*>
    \\&\hspace{0.4cm} - 2\gamma\< \bar F^{k+1/2}_{m_k} - F(\hat u^k_{m_k}) - F_{m_k} (\tilde u^k_{m_k}) + F_{m_k} (\hat u^k_{m_k}), \hat u^{k}_{m_k} - z^*> \\&\hspace{0.4cm}+ \sqn{\tilde z^{k+1}_{m_k} - \hat u^{k}_{m_k}} - \sqn{z^{k}_{m_k} - \hat u^{k}_{m_k}}
	\\&\hspace{0.4cm} + 2\<\tilde u^k_{m_k} - \hat u^{k}_{m_k}, \hat u^{k}_{m_k} - z^*>  + 2\< \hat z^{k+1}_{m_k} - \tilde z^{k+1}_{m_k}, \tilde z^{k+1}_{m_k} - z^* > + \sqn{\tilde z^{k+1}_{m_k} - \hat z^{k+1}_{m_k}}.
	\end{align*}
Applying property of the solution $z^*$: $\<\gamma F(z^*), z^* - z > \leq 0$ (for all $z \in \mathcal{Z}$). And then $\mu$-strong convexity - strong concavity of $f$, we obtain
    \begin{align*}
	\sqn{z^{k+1}_{m_k} - z^*}
	&\leq
	\sqn{z^{k}_{m_k} - z^*} - 2\gamma\<F(\hat u^{k}_{m_k}) - F(z^*), \hat u^{k}_{m_k} - z^*> 
    \\&\hspace{0.4cm} - 2\gamma\< \bar F^{k+1/2}_{m_k} - F(\hat u^k_{m_k}) - F_{m_k} (\tilde u^k_{m_k}) + F_{m_k} (\hat u^k_{m_k}), \hat u^{k}_{m_k} - z^*> 
    \\&\hspace{0.4cm} + \sqn{\tilde z^{k+1}_{m_k} - \hat u^{k}_{m_k}} - \sqn{z^{k}_{m_k} - \hat u^{k}_{m_k}}
	\\&\hspace{0.4cm} + 2\<\tilde u^k_{m_k} - \hat u^{k}_{m_k}, \hat u^{k}_{m_k} - z^*>  + 2\< \hat z^{k+1}_{m_k} - \tilde z^{k+1}_{m_k}, \tilde z^{k+1}_{m_k} - z^* > + \sqn{\tilde z^{k+1}_{m_k} - \hat z^{k+1}_{m_k}}
	\\&\leq
	\sqn{z^{k}_{m_k} - z^*} - 2\gamma \mu \sqn{ \hat u^{k}_{m_k} - z^*}
    \\&\hspace{0.4cm} - 2\gamma\< \bar F^{k+1/2}_{m_k} - F(\hat u^k_{m_k}) - F_{m_k} (\tilde u^k_{m_k}) + F_{m_k} (\hat u^k_{m_k}), \hat u^{k}_{m_k} - z^*> 
    \\&\hspace{0.4cm} + \sqn{\tilde z^{k+1}_{m_k} - \hat u^{k}_{m_k}} - \sqn{z^{k}_{m_k} - \hat u^{k}_{m_k}}
	\\&\hspace{0.4cm} + 2\<\tilde u^k_{m_k} - \hat u^{k}_{m_k}, \hat u^{k}_{m_k} - z^*>  + 2\< \hat z^{k+1}_{m_k} - \tilde z^{k+1}_{m_k}, \tilde z^{k+1}_{m_k} - z^* > + \sqn{\tilde z^{k+1}_{m_k} - \hat z^{k+1}_{m_k}}.
	\end{align*}
By Young's inequality, we have
	\begin{align*}
	\sqn{z^{k+1}_{m_k} - z^*}
	&\leq
	\sqn{z^{k}_{m_k} - z^*} - 2\gamma \mu \sqn{ \hat u^{k}_{m_k} - z^*}
    \\&\hspace{0.4cm} + \frac{4\gamma}{\mu}\sqn{\bar F^{k+1/2}_{m_k} - F(\hat u^k_{m_k}) - F_{m_k} (\tilde u^k_{m_k}) + F_{m_k} (\hat u^k_{m_k})}  
    \\&\hspace{0.4cm}+ \frac{\gamma \mu}{4}\sqn{\hat u^{k}_{m_k} - z^*} + \sqn{\tilde z^{k+1}_{m_k} - \hat u^{k}_{m_k}} - \sqn{z^{k}_{m_k} - \hat u^{k}_{m_k}}
	\\&\hspace{0.4cm} + \frac{4}{\gamma\mu}\sqn{\tilde u^k_{m_k} - \hat u^{k}_{m_k}}  + \frac{\gamma \mu}{4}\sqn{\hat u^{k}_{m_k} - z^*} 
	\\&\hspace{0.4cm} + 2\< \hat z^{k+1}_{m_k} - \tilde z^{k+1}_{m_k}, \tilde z^{k+1}_{m_k} - z^* > + \sqn{\tilde z^{k+1}_{m_k} - \hat z^{k+1}_{m_k}} \\&= 
	\sqn{z^{k}_{m_k} - z^*} - \frac{3\gamma \mu}{2}\sqn{ \hat u^{k}_{m_k} - z^*}
    \\&\hspace{0.4cm} + \frac{4\gamma}{\mu}\sqn{\bar F^{k+1/2}_{m_k} - F(\hat u^k_{m_k}) - F_{m_k} (\tilde u^k_{m_k}) + F_{m_k} (\hat u^k_{m_k})} 
    \\&\hspace{0.4cm} + \sqn{\tilde u^k_{m_k} + \gamma \cdot (\bar F^k_{m_k} - F_{m_k} (z^k_{m_k}) - \bar F^{k+1/2}_{m_k} + F_{m_k} (\tilde u^k_{m_k})) - \hat u^{k}_{m_k}} 
    \\&\hspace{0.4cm} - \sqn{z^{k}_{m_k} - \hat u^{k}_{m_k}}
	+ \frac{4}{\gamma \mu }\sqn{\tilde u^k_{m_k} - \hat u^{k}_{m_k}} 
	\\&\hspace{0.4cm} + 2\< \hat z^{k+1}_{m_k} - \tilde z^{k+1}_{m_k}, \tilde z^{k+1}_{m_k} - z^* > + \sqn{\tilde z^{k+1}_{m_k} - \hat z^{k+1}_{m_k}}
	\\&\leq 
	\sqn{z^{k}_{m_k} - z^*} - \frac{3\gamma \mu}{2} \sqn{ \hat u^{k}_{m_k} - z^*}
    \\&\hspace{0.4cm} + \frac{8\gamma}{\mu}\sqn{ F(\tilde u^k_{m_k}) - F(\hat u^k_{m_k}) - F_{m_k} (\tilde u^k_{m_k}) + F_{m_k} (\hat u^k_{m_k})} 
    \\&\hspace{0.4cm} + \frac{8\gamma}{\mu}\sqn{\bar F^{k+1/2}_{m_k} - F(\tilde u^k_{m_k})} 
    \\&\hspace{0.4cm} + 6\gamma^2\sqn{ F(z^k_{m_k}) - F_{m_k} (z^k_{m_k}) - F(\tilde u^k_{m_k}) + F_{m_k} (\tilde u^k_{m_k})} 
    \\&\hspace{0.4cm} + 6\gamma^2\sqn{ \bar F^k_{m_k} - F(z^k_{m_k})} + 6\gamma^2\sqn{ \bar F^{k+1/2}_{m_k} - F(\tilde u^k_{m_k})} 
    \\&\hspace{0.4cm} + 2\sqn{\tilde u^k_{m_k} - \hat u^{k}_{m_k}} - \sqn{z^{k}_{m_k} - \hat u^{k}_{m_k}}
	\\&\hspace{0.4cm} + \frac{4}{\gamma \mu }\sqn{\tilde u^k_{m_k} - \hat u^{k}_{m_k}} + 2\< \hat z^{k+1}_{m_k} - \tilde z^{k+1}_{m_k}, \tilde z^{k+1}_{m_k} - z^* > 
	+ \sqn{\tilde z^{k+1}_{m_k} - \hat z^{k+1}_{m_k}}.
	\end{align*}
Note that the function $f-f_{m_k}$ is $\delta$ - smooth (since $\sqn{\nabla_{xx} f- \nabla_{xx}f_{m_k}} \leq \delta$, $\sqn{\nabla_{xy} f- \nabla_{xy}f_{m_k}} \leq \delta$, $\sqn{\nabla_{yy} f- \nabla_{yy}f_{m_k}} \leq \delta$), then
    \begin{align*}
	\sqn{z^{k+1}_{m_k} - z^*}
	&\leq 
	\sqn{z^{k}_{m_k} - z^*} - \frac{3\gamma \mu}{2} \sqn{ \hat u^{k}_{m_k} - z^*}
    \\&\hspace{0.4cm} + \frac{8\gamma \delta^2}{\mu}\sqn{ \tilde u^k_{m_k}- \hat u^k_{m_k}} + \frac{8\gamma}{\mu}\sqn{\bar F^{k+1/2}_{m_k} - F(\tilde u^k_{m_k})} 
    \\&\hspace{0.4cm} + 6\gamma^2\delta^2 \sqn{ z^k_{m_k}-\tilde u^k_{m_k}} + 6\gamma^2\sqn{ \bar F^k_{m_k} - F(z^k_{m_k})} + 6\gamma^2\sqn{ \bar F^{k+1/2}_{m_k} - F(\tilde u^k_{m_k})}
    \\&\hspace{0.4cm}  + 2\sqn{\tilde u^k_{m_k} - \hat u^{k}_{m_k}} - \sqn{z^{k}_{m_k} - \hat u^{k}_{m_k}}
	\\&\hspace{0.4cm} + \frac{4}{\gamma \mu }\sqn{\tilde u^k_{m_k} - \hat u^{k}_{m_k}} + 2\< \hat z^{k+1}_{m_k} - \tilde z^{k+1}_{m_k}, \tilde z^{k+1}_{m_k} - z^* > + \sqn{\tilde z^{k+1}_{m_k} - \hat z^{k+1}_{m_k}}
	\\&\leq
	\sqn{z^{k}_{m_k} - z^*} - \frac{3\gamma \mu}{2} \sqn{ \hat u^{k}_{m_k} - z^*} - (1 - 12\gamma^2\delta^2)\sqn{z^{k}_{m_k} - \hat u^{k}_{m_k}}
    \\&\hspace{0.4cm} + \left(2 + 12\gamma^2\delta^2 +\frac{4}{\gamma \mu }+ \frac{8\gamma \delta^2}{\mu}\right)\sqn{ \tilde u^k_{m_k}- \hat u^k_{m_k}} 
    \\&\hspace{0.4cm} + 6\gamma^2\sqn{ \bar F^k_{m_k} - F(z^k_{m_k})} + \left(6\gamma^2 + \frac{8\gamma}{\mu}\right)\sqn{ \bar F^{k+1/2}_{m_k} - F(\tilde u^k_{m_k})}
	\\&\hspace{0.4cm} + 2\< \hat z^{k+1}_{m_k} - \tilde z^{k+1}_{m_k}, \tilde z^{k+1}_{m_k} - z^* > + \sqn{\tilde z^{k+1}_{m_k} - \hat z^{k+1}_{m_k}}.
	\end{align*}
	By inequality $\sqn{a+b} \geq \frac{2}{3}\sqn{a} - 2\sqn{b}$, we have
	\begin{align*}
	\sqn{z^{k+1}_{m_k} - z^*}
	&\leq
	(1 - \gamma \mu)\sqn{z^{k}_{m_k} - z^*} - (1 - 3\gamma \mu - 12\gamma^2\delta^2)\sqn{z^{k}_{m_k} - \hat u^{k}_{m_k}}
    \\&\hspace{0.4cm} + \left(2 + 12\gamma^2\delta^2+ \frac{4}{\gamma \mu } + \frac{8\gamma \delta^2}{\mu}\right)\sqn{ \tilde u^k_{m_k}- \hat u^k_{m_k}} 
    \\&\hspace{0.4cm} + 6\gamma^2\sqn{ \bar F^k_{m_k} - F(z^k_{m_k})} + \left(6\gamma^2 + \frac{8\gamma}{\mu}\right)\sqn{ \bar F^{k+1/2}_{m_k} - F(\tilde u^k_{m_k})}
	\\&\hspace{0.4cm} + 2\< \hat z^{k+1}_{m_k} - \tilde z^{k+1}_{m_k}, \tilde z^{k+1}_{m_k} - z^* > + \sqn{\tilde z^{k+1}_{m_k} - \hat z^{k+1}_{m_k}}.
	\end{align*}
\EndProof

\begin{lemma} \label{lem5}
	Let for problem \eqref{inner} we use Extragradient method with starting point $z^{k}_{m_k}$ and number of iterations:
	\begin{equation}
	\label{T3}
		T = \cO  \left(\left(1 + \gamma L\right)\log\frac{1}{\tilde e} \right).
	\end{equation}
	Then for an output $\tilde u^k_{m_k}$ it holds that 
	\begin{equation*}
		\sqn{\tilde u^k_{m_k} - \hat{u}^k_{m_k}}\leq \tilde e\sqn{z^{k}_{m_k} - \hat u^k_{m_k}}.
	\end{equation*}
\end{lemma}

\begin{theorem}
Let problem \eqref{inner} be solved by Extragradient with precision $\tilde e$:
\begin{align}
\label{e3}
    \tilde e = \frac{1}{2\left(2 + 12\gamma^2\delta^2+ \frac{4}{\gamma \mu } + \frac{8\gamma \delta^2}{\mu}\right)}
\end{align}
and number of iterations $T$ from \eqref{T3}. Suppose that parameters $H_0$ and $H_1$ satisfy
\begin{align}
\label{H}
    H_0 = \mathcal{O} \left(\frac{1}{\sqrt{\rho}}\log \left(\frac{\left(\gamma^2 + \frac{\gamma}{\mu}\right) \cdot M(L\Omega + G)^2}{\varepsilon \gamma \mu}\right)\right), \notag\\ 
    H_1 = \mathcal{O} \left(\frac{1}{\sqrt{\rho}}\log \left(\frac{\left(1 + \gamma^2 L^2 + \frac{\gamma L^2}{\mu}\right) \cdot M\Omega^2}{\varepsilon \gamma \mu}\right)\right)
\end{align}
Additionally, let us choose stepsize $\gamma$ as follows
\begin{align}
\label{gamma3}
    \gamma = \min\left\{ \frac{1}{12 \mu}; \frac{1}{7 \delta}\right\}.
\end{align}
Then Algorithm 1 converges linearly to the solution $z^*$ and it holds that $\sqn{ z^{K} - z^*} \sim \varepsilon$ after 
\begin{align}
\label{K3}
    K = \cO \left(\frac{1}{\gamma \mu} \log \frac{\sqn{z^{0} - z^*}}{\varepsilon} \right)\quad \text{iterations}.
\end{align}
\end{theorem}

\textbf{Proof:} Combining results from Lemma \ref{lem4} and \ref{lem5} gives
\begin{align*}
	\sqn{z^{k+1}_{m_k} - z^*}
	&\leq
	(1 - \gamma \mu)\sqn{z^{k}_{m_k} - z^*} - (1 - 3\gamma \mu - 12\gamma^2\delta^2)\sqn{z^{k}_{m_k} - \hat u^{k}_{m_k}}
    \\&\hspace{0.4cm} + \left(2 + 12\gamma^2\delta^2+ \frac{4}{\gamma \mu } + \frac{8\gamma \delta^2}{\mu}\right)\sqn{ \tilde u^k_{m_k}- \hat u^k_{m_k}} 
    \\&\hspace{0.4cm} + 6\gamma^2\sqn{ \bar F^k_{m_k} - F(z^k_{m_k})} + \left(6\gamma^2 + \frac{8\gamma}{\mu}\right)\sqn{ \bar F^{k+1/2}_{m_k} - F(\tilde u^k_{m_k})}
	\\&\hspace{0.4cm} + 2\< \hat z^{k+1}_{m_k} - \tilde z^{k+1}_{m_k}, \tilde z^{k+1}_{m_k} - z^* > + \sqn{\tilde z^{k+1}_{m_k} - \hat z^{k+1}_{m_k}}.
\end{align*}
With the choice $e$ from \eqref{e3} and $\gamma$ from \eqref{gamma3}, we obtain
\begin{align*}
	\sqn{z^{k+1}_{m_k} - z^*}
	&\leq
	(1 - \gamma \mu)\sqn{z^{k}_{m_k} - z^*}
	\\&\hspace{0.4cm} + 6\gamma^2\sqn{ \bar F^k_{m_k} - F(z^k_{m_k})} + \left(6\gamma^2 + \frac{8\gamma}{\mu}\right)\sqn{ \bar F^{k+1/2}_{m_k} - F(\tilde u^k_{m_k})}
	\\&\hspace{0.4cm} + 2\< \hat z^{k+1}_{m_k} - \tilde z^{k+1}_{m_k}, \tilde z^{k+1}_{m_k} - z^* > + \sqn{\tilde z^{k+1}_{m_k} - \hat z^{k+1}_{m_k}}.
\end{align*}
Passing from the local $z^{k+1}_{m_k}$ and $z^{k}_{m_k}$ to $\bar z^{k+1}$ and $\bar z^{k}$, we have
\begin{align}
    \label{505}
	\sqn{\bar z^{k+1} - z^*}
	&\leq \notag
	(1 - \gamma \mu)\sqn{\bar z^{k} - z^*}
	\\&\hspace{0.4cm} + 6\gamma^2\sqn{ \bar F^k_{m_k} - F(z^k_{m_k})} + \left(6\gamma^2 + \frac{8\gamma}{\mu}\right)\sqn{ \bar F^{k+1/2}_{m_k} - F(\tilde u^k_{m_k})} \notag
	\\&\hspace{0.4cm} + 2\| \hat z^{k+1}_{m_k} - \tilde z^{k+1}_{m_k}\| \cdot \| \tilde z^{k+1}_{m_k} - z^* \| + \sqn{\tilde z^{k+1}_{m_k} - \hat z^{k+1}_{m_k}} \notag
	\\&\hspace{0.4cm} + 2 \| z^{k+1}_{m_k} - \bar z^{k+1} \| \cdot  \| \bar z^{k+1} - z^* \| +\sqn{z^{k+1}_{m_k} - \bar z^{k+1}} \notag
	\\&\hspace{0.4cm} + 2\|z^{k}_{m_k} - \bar z^{k} \| \cdot  \| \bar z^{k} - z^*\| +\sqn{z^{k}_{m_k} - \bar z^{k}}.
\end{align}
Further we will work separately only with the last 4 lines, because the last 4 lines depend on the number of iterations $H_0$ and $H_1$, then we can make them small by choosing the correct $H_0$ and $H_1$.
\begin{align*}
	\text{Err}(k)
	&=
	6\gamma^2\sqn{ \bar F^k_{m_k} - F(z^k_{m_k})} + \left(6\gamma^2 + \frac{8\gamma}{\mu}\right)\sqn{ \bar F^{k+1/2}_{m_k} - F(\tilde u^k_{m_k})}
	\\&\hspace{0.4cm} + 2\| \hat z^{k+1}_{m_k} - \tilde z^{k+1}_{m_k}\| \cdot \| \tilde z^{k+1}_{m_k} - z^* \| + \sqn{\tilde z^{k+1}_{m_k} - \hat z^{k+1}_{m_k}}
	\\&\hspace{0.4cm} + 2 \| z^{k+1}_{m_k} - \bar z^{k+1} \| \cdot  \| \bar z^{k+1} - z^* \| +\sqn{z^{k+1}_{m_k} - \bar z^{k+1}}
	\\&\hspace{0.4cm} + 2\|z^{k}_{m_k} - \bar z^{k} \| \cdot  \| \bar z^{k} - z^*\| +\sqn{z^{k}_{m_k} - \bar z^{k}}
	\\&\leq
	6\gamma^2\sqn{ \bar F^k_{m_k} - F(z^k_{m_k})} + \left(6\gamma^2 + \frac{8\gamma}{\mu}\right)\sqn{ \bar F^{k+1/2}_{m_k} - F(\tilde u^k_{m_k})}
	\\&\hspace{0.4cm} + 2\| \hat z^{k+1}_{m_k} - \tilde z^{k+1}_{m_k}\| \cdot \Omega + \sqn{\tilde z^{k+1}_{m_k} - \hat z^{k+1}_{m_k}}
	\\&\hspace{0.4cm} + 2 \| z^{k+1}_{m_k} - \bar z^{k+1} \| \cdot  \Omega +\sqn{z^{k+1}_{m_k} - \bar z^{k+1}}
	\\&\hspace{0.4cm} + 2\|z^{k}_{m_k} - \bar z^{k} \| \cdot  \Omega +\sqn{z^{k}_{m_k} - \bar z^{k}}
\end{align*}	
Next we use the definition of $\bar z^{k}$ and $\bar z^{k+1}$ and the fact from line 6 of Algorithm 2: $M \tilde z^{k+1}_{m_k} = \sum_{i=1}^M \hat z^{k+1}_{i}$, and get
\begin{align*}
	\text{Err}(k)
    &\leq
	12\gamma^2\sqn{ \bar F^k_{m_k} - \frac{1}{M} \sum\limits_{i=1}^M F_i(z^k_i)} + 12\gamma^2\sqn{\frac{1}{M} \sum\limits_{i=1}^M F_i(z^k_i)  - F(z^k_{m_k})} 
	\\&\hspace{0.4cm}+ \left(12\gamma^2 + \frac{16\gamma}{\mu}\right)\sqn{ \bar F^{k+1/2}_{m_k} - \frac{1}{M} \sum\limits_{i=1}^M F_i(u^k_i)} 
	\\&\hspace{0.4cm}+ \left(12\gamma^2 + \frac{16\gamma}{\mu}\right)\sqn{ \frac{1}{M} \sum\limits_{i=1}^M F_i(u^k_i) - F(\tilde u^k_{m_k})}
	\\&\hspace{0.4cm} + 2\left\| \hat z^{k+1}_{m_k} - \frac{1}{M} \sum\limits_{i=1}^M\hat z^{k+1}_{i}\right\| \cdot \Omega + \sqn{\frac{1}{M} \sum\limits_{i=1}^M\hat z^{k+1}_{i} - \hat z^{k+1}_{m_k}}
	\\&\hspace{0.4cm} + 2 \left\| \text{proj}[\hat z^{k+1}_{m_k}] - \frac{1}{M} \sum\limits_{i=1}^M \text{proj}[\hat z^{k+1}_{i}] \right\| \cdot  \Omega +\sqn{\text{proj}[\hat z^{k+1}_{m_k}]  -\frac{1}{M} \sum\limits_{i=1}^M \text{proj}[\hat z^{k+1}_{i}]}
	\\&\hspace{0.4cm} + 2\left\| \text{proj}[\hat z^{k}_{m_k}] - \frac{1}{M} \sum\limits_{i=1}^M \text{proj}[\hat z^{k}_{i}] \right\| \cdot  \Omega +\sqn{\text{proj}[\hat z^{k}_{m_k}] - \frac{1}{M} \sum\limits_{i=1}^M \text{proj}[\hat z^{k}_{i}] }
	\\&\leq
	12\gamma^2\sqn{ \bar F^k_{m_k} - \frac{1}{M} \sum\limits_{i=1}^M F_i(z^k_i)} + 12\gamma^2\sqn{\frac{1}{M} \sum\limits_{i=1}^M F_i(z^k_i)  - F(z^k_{m_k})} 
	\\&\hspace{0.4cm}+ \left(12\gamma^2 + \frac{16\gamma}{\mu}\right)\sqn{ \bar F^{k+1/2}_{m_k} - \frac{1}{M} \sum\limits_{i=1}^M F_i(u^k_i)} 
	\\&\hspace{0.4cm}+ \left(12\gamma^2 + \frac{16\gamma}{\mu}\right)\sqn{ \frac{1}{M} \sum\limits_{i=1}^M F_i(u^k_i) - F(\tilde u^k_{m_k})}
	\\&\hspace{0.4cm} + 2\left\| \hat z^{k+1}_{m_k} - \frac{1}{M} \sum\limits_{i=1}^M\hat z^{k+1}_{i}\right\| \cdot \Omega + \sqn{\frac{1}{M} \sum\limits_{i=1}^M\hat z^{k+1}_{i} - \hat z^{k+1}_{m_k}}
	\\&\hspace{0.4cm} + \frac{2}{M} \sum\limits_{i=1}^M \left\| \hat z^{k+1}_{m_k} - \hat z^{k+1}_{i} \right\| \cdot  \Omega +\frac{1}{M} \sum\limits_{i=1}^M \sqn{\hat z^{k+1}_{m_k} - \hat z^{k+1}_{i}}
	\\&\hspace{0.4cm} + \frac{2}{M} \sum\limits_{i=1}^M \left\| \hat z^{k}_{m_k} -  \hat z^{k}_{i} \right\| \cdot  \Omega +\frac{1}{M} \sum\limits_{i=1}^M \sqn{\hat z^{k}_{m_k} -  \hat z^{k}_{i}}\\
	&\leq
	12\gamma^2\sqn{ \bar F^k_{m_k} - \frac{1}{M} \sum\limits_{i=1}^M F_i(z^k_i)} + \left(12\gamma^2 + \frac{16\gamma}{\mu}\right)\sqn{ \bar F^{k+1/2}_{m_k} - \frac{1}{M} \sum\limits_{i=1}^M F_i(u^k_i)}
	\\&\hspace{0.4cm} + 2\Omega\left\| \hat z^{k+1}_{m_k} - \frac{1}{M} \sum\limits_{i=1}^M\hat z^{k+1}_{i}\right\| + \sqn{\frac{1}{M} \sum\limits_{i=1}^M\hat z^{k+1}_{i} - \hat z^{k+1}_{m_k}}
	\\&\hspace{0.4cm} + 12\gamma^2 \frac{1}{M} \sum\limits_{i=1}^M \sqn{ F_i(z^k_i)  - F_i(z^k_{m_k})}  + \left(12\gamma^2 + \frac{16\gamma}{\mu}\right) \frac{1}{M} \sum\limits_{i=1}^M \sqn{  F_i(u^k_i) - F_i(\tilde u^k_{m_k})}
	\\&\hspace{0.4cm} + \frac{2}{M} \sum\limits_{i=1}^M \left\| \hat z^{k+1}_{m_k} - \frac{1}{M} \sum\limits_{j=1}^M \hat z^{k+1}_{j} + \frac{1}{M} \sum\limits_{j=1}^M \hat z^{k+1}_{j} - \hat z^{k+1}_{i} \right\| \cdot  \Omega 
	\\&\hspace{0.4cm}+\frac{1}{M} \sum\limits_{i=1}^M \sqn{\hat z^{k+1}_{m_k} - \frac{1}{M} \sum\limits_{j=1}^M \hat z^{k+1}_{j} + \frac{1}{M} \sum\limits_{j=1}^M \hat z^{k+1}_{j} - \hat z^{k+1}_{i}}
\end{align*}
\begin{align*}
	\\&\hspace{0.4cm} + \frac{2}{M} \sum\limits_{i=1}^M \left\| \hat z^{k}_{m_k} - \frac{1}{M} \sum\limits_{j=1}^M \hat z^{k}_{j} + \frac{1}{M} \sum\limits_{j=1}^M \hat z^{k}_{j} -  \hat z^{k}_{i} \right\| \cdot  \Omega 	\\&\hspace{0.4cm} +\frac{1}{M} \sum\limits_{i=1}^M \sqn{\hat z^{k}_{m_k} - \frac{1}{M} \sum\limits_{j=1}^M \hat z^{k}_{j} + \frac{1}{M} \sum\limits_{j=1}^M \hat z^{k}_{j} -  \hat z^{k}_{i}}
	\\&\leq
	12\gamma^2\sqn{ \bar F^k_{m_k} - \frac{1}{M} \sum\limits_{i=1}^M F_i(z^k_i)} + \left(12\gamma^2 + \frac{16\gamma}{\mu}\right)\sqn{ \bar F^{k+1/2}_{m_k} - \frac{1}{M} \sum\limits_{i=1}^M F_i(u^k_i)}
	\\&\hspace{0.4cm} + 2\Omega\left\| \hat z^{k+1}_{m_k} - \frac{1}{M} \sum\limits_{i=1}^M\hat z^{k+1}_{i}\right\| + \sqn{\frac{1}{M} \sum\limits_{i=1}^M\hat z^{k+1}_{i} - \hat z^{k+1}_{m_k}}
	\\&\hspace{0.4cm} + 12\gamma^2 L^2 \frac{1}{M} \sum\limits_{i=1}^M \sqn{ z^k_i  - z^k_{m_k}}  + \left(12\gamma^2 + \frac{16\gamma}{\mu}\right) L^2 \frac{1}{M} \sum\limits_{i=1}^M \sqn{  u^k_i - \tilde u^k_{m_k}}
	\\&\hspace{0.4cm} + \frac{2}{M} \sum\limits_{i=1}^M \left\| \hat z^{k+1}_{m_k} - \frac{1}{M} \sum\limits_{j=1}^M \hat z^{k+1}_{j} \right\| \cdot  \Omega + \frac{2}{M} \sum\limits_{i=1}^M \left\| \frac{1}{M} \sum\limits_{j=1}^M \hat z^{k+1}_{j} - \hat z^{k+1}_{i} \right\| \cdot  \Omega
	\\&\hspace{0.4cm}+\frac{2}{M} \sum\limits_{i=1}^M \sqn{\hat z^{k+1}_{m_k} - \frac{1}{M} \sum\limits_{j=1}^M \hat z^{k+1}_{j}} + \frac{2}{M} \sum\limits_{i=1}^M \sqn{\frac{1}{M} \sum\limits_{j=1}^M \hat z^{k+1}_{j} - \hat z^{k+1}_{i}}
	\\&\hspace{0.4cm} + \frac{2}{M} \sum\limits_{i=1}^M \left\| \hat z^{k}_{m_k} - \frac{1}{M} \sum\limits_{j=1}^M \hat z^{k}_{j} \right\| \cdot  \Omega + \frac{2}{M} \sum\limits_{i=1}^M \left\| \frac{1}{M} \sum\limits_{j=1}^M \hat z^{k}_{j} -  \hat z^{k}_{i} \right\| \cdot  \Omega 
	\\&\hspace{0.4cm} +\frac{2}{M} \sum\limits_{i=1}^M \sqn{\hat z^{k}_{m_k} - \frac{1}{M} \sum\limits_{j=1}^M \hat z^{k}_{j}} +\frac{2}{M} \sum\limits_{i=1}^M \sqn{\frac{1}{M} \sum\limits_{j=1}^M \hat z^{k}_{j} -  \hat z^{k}_{i}}.
\end{align*}
Small rearrangement gives
\begin{align*}
	\text{Err}(k)
	&\leq
	12\gamma^2\sqn{ \bar F^k_{m_k} - \frac{1}{M} \sum\limits_{i=1}^M F_i(z^k_i)} + \left(12\gamma^2 + \frac{16\gamma}{\mu}\right)\sqn{ \bar F^{k+1/2}_{m_k} - \frac{1}{M} \sum\limits_{i=1}^M F_i(u^k_i)}
	\\&\hspace{0.4cm} + 4\Omega\left\| \hat z^{k+1}_{m_k} - \frac{1}{M} \sum\limits_{i=1}^M\hat z^{k+1}_{i}\right\| + 4\sqn{\frac{1}{M} \sum\limits_{i=1}^M\hat z^{k+1}_{i} - \hat z^{k+1}_{m_k}} + 2 \Omega \left\| \hat z^{k}_{m_k} - \frac{1}{M} \sum\limits_{j=1}^M \hat z^{k}_{j} \right\|
	\\&\hspace{0.4cm}  + 2 \sqn{\hat z^{k}_{m_k} - \frac{1}{M} \sum\limits_{j=1}^M \hat z^{k}_{j}}
	+ 12\gamma^2 L^2 \frac{1}{M} \sum\limits_{i=1}^M \sqn{ \text{proj}[\hat z^k_i] -  \text{proj}[\hat z^k_{m_k}]}  
	\\&\hspace{0.4cm} + \left(12\gamma^2 + \frac{16\gamma}{\mu}\right) L^2 \frac{1}{M} \sum\limits_{i=1}^M \sqn{  u^k_i - \frac{1}{M} \sum\limits_{j=1}^M u^k_j + \frac{1}{M} \sum\limits_{j=1}^M u^k_j - \tilde u^k_{m_k}}
	\\&\hspace{0.4cm} + \frac{2}{M} \sum\limits_{i=1}^M \left\| \frac{1}{M} \sum\limits_{j=1}^M \hat z^{k+1}_{j} - \hat z^{k+1}_{i} \right\| \cdot  \Omega
	+ \frac{2}{M} \sum\limits_{i=1}^M \sqn{\frac{1}{M} \sum\limits_{j=1}^M \hat z^{k+1}_{j} - \hat z^{k+1}_{i}}
	\\&\hspace{0.4cm} + \frac{2}{M} \sum\limits_{i=1}^M \left\| \frac{1}{M} \sum\limits_{j=1}^M \hat z^{k}_{j} -  \hat z^{k}_{i} \right\| \cdot  \Omega  +\frac{2}{M} \sum\limits_{i=1}^M \sqn{\frac{1}{M} \sum\limits_{j=1}^M \hat z^{k}_{j} -  \hat z^{k}_{i}}.
\end{align*}
\begin{align*}
    \text{Err}(k)	
	&\leq
	12\gamma^2\sqn{ \bar F^k_{m_k} - \frac{1}{M} \sum\limits_{i=1}^M F_i(z^k_i)} + \left(12\gamma^2 + \frac{16\gamma}{\mu}\right)\sqn{ \bar F^{k+1/2}_{m_k} - \frac{1}{M} \sum\limits_{i=1}^M F_i(u^k_i)}
	\\&\hspace{0.4cm} + 4\Omega\left\| \hat z^{k+1}_{m_k} - \frac{1}{M} \sum\limits_{i=1}^M\hat z^{k+1}_{i}\right\| + 4\sqn{\frac{1}{M} \sum\limits_{i=1}^M\hat z^{k+1}_{i} - \hat z^{k+1}_{m_k}}
	\\&\hspace{0.4cm} + 2 \Omega \left\| \hat z^{k}_{m_k} - \frac{1}{M} \sum\limits_{j=1}^M \hat z^{k}_{j} \right\| + 2 \sqn{\hat z^{k}_{m_k} - \frac{1}{M} \sum\limits_{j=1}^M \hat z^{k}_{j}}
	\\&\hspace{0.4cm} + 12\gamma^2 L^2 \frac{1}{M} \sum\limits_{i=1}^M \sqn{ \hat z^k_i -\frac{1}{M} \sum\limits_{j=1}^M \hat z^{k}_{j} + \frac{1}{M} \sum\limits_{j=1}^M \hat z^{k}_{j} -  \hat z^k_{m_k}}  
	\\&\hspace{0.4cm} + \left(24\gamma^2 + \frac{32\gamma}{\mu}\right) L^2 \frac{1}{M} \sum\limits_{i=1}^M \sqn{  u^k_i - \frac{1}{M} \sum\limits_{j=1}^M u^k_j} 
	\\&\hspace{0.4cm}+ \left(24\gamma^2 + \frac{32\gamma}{\mu}\right) L^2 \sqn{  \frac{1}{M} \sum\limits_{j=1}^M u^k_j - \tilde u^k_{m_k}}
	\\&\hspace{0.4cm} + \frac{2}{M} \sum\limits_{i=1}^M \left\| \frac{1}{M} \sum\limits_{j=1}^M \hat z^{k+1}_{j} - \hat z^{k+1}_{i} \right\| \cdot  \Omega
	+ \frac{2}{M} \sum\limits_{i=1}^M \sqn{\frac{1}{M} \sum\limits_{j=1}^M \hat z^{k+1}_{j} - \hat z^{k+1}_{i}}
	\\&\hspace{0.4cm} + \frac{2}{M} \sum\limits_{i=1}^M \left\| \frac{1}{M} \sum\limits_{j=1}^M \hat z^{k}_{j} -  \hat z^{k}_{i} \right\| \cdot  \Omega  +\frac{2}{M} \sum\limits_{i=1}^M \sqn{\frac{1}{M} \sum\limits_{j=1}^M \hat z^{k}_{j} -  \hat z^{k}_{i}}
	\\&\leq
	12\gamma^2\sqn{ \bar F^k_{m_k} - \frac{1}{M} \sum\limits_{i=1}^M F_i(z^k_i)} + \left(12\gamma^2 + \frac{16\gamma}{\mu}\right)\sqn{ \bar F^{k+1/2}_{m_k} - \frac{1}{M} \sum\limits_{i=1}^M F_i(u^k_i)}
	\\&\hspace{0.4cm} + 4\Omega\left\| \hat z^{k+1}_{m_k} - \frac{1}{M} \sum\limits_{i=1}^M\hat z^{k+1}_{i}\right\| + 4\sqn{\frac{1}{M} \sum\limits_{i=1}^M\hat z^{k+1}_{i} - \hat z^{k+1}_{m_k}}
	\\&\hspace{0.4cm} + 2 \Omega \left\| \hat z^{k}_{m_k} - \frac{1}{M} \sum\limits_{j=1}^M \hat z^{k}_{j} \right\| + (2 + 24\gamma^2 L^2) \sqn{\hat z^{k}_{m_k} - \frac{1}{M} \sum\limits_{j=1}^M \hat z^{k}_{j}} 
	\\&\hspace{0.4cm}+ 24\gamma^2 L^2 \frac{1}{M} \sum\limits_{i=1}^M \sqn{ \hat z^k_i -\frac{1}{M} \sum\limits_{j=1}^M \hat z^{k}_{j}}   
	\\&\hspace{0.4cm} + \left(24\gamma^2 + \frac{32\gamma}{\mu}\right) L^2 \frac{1}{M} \sum\limits_{i=1}^M \sqn{  u^k_i - \frac{1}{M} \sum\limits_{j=1}^M u^k_j} 
	\\&\hspace{0.4cm}+ \left(24\gamma^2 + \frac{32\gamma}{\mu}\right) L^2 \sqn{  \frac{1}{M} \sum\limits_{j=1}^M u^k_j - \tilde u^k_{m_k}}
	\\&\hspace{0.4cm} + \frac{2}{M} \sum\limits_{i=1}^M \left\| \frac{1}{M} \sum\limits_{j=1}^M \hat z^{k+1}_{j} - \hat z^{k+1}_{i} \right\| \cdot  \Omega
	+ \frac{2}{M} \sum\limits_{i=1}^M \sqn{\frac{1}{M} \sum\limits_{j=1}^M \hat z^{k+1}_{j} - \hat z^{k+1}_{i}}
	\\&\hspace{0.4cm} + \frac{2}{M} \sum\limits_{i=1}^M \left\| \frac{1}{M} \sum\limits_{j=1}^M \hat z^{k}_{j} -  \hat z^{k}_{i} \right\| \cdot  \Omega  +\frac{2}{M} \sum\limits_{i=1}^M \sqn{\frac{1}{M} \sum\limits_{j=1}^M \hat z^{k}_{j} -  \hat z^{k}_{i}}.
\end{align*}
Now we are ready to apply {\tt AccGossip} convergence results (\eqref{400}, \eqref{401}, \eqref{402}) to each of these terms:
\begin{align*}
	\text{Err}(k)
	&\leq
	12\gamma^2 (1 - \sqrt{\rho})^{2H_0} \cdot 2M(L\Omega + G)^2 + \left(12\gamma^2 + \frac{16\gamma}{\mu}\right)(1 - \sqrt{\rho})^{2H_0} \cdot 2M(L\Omega + G)^2
	\\&\hspace{0.4cm} + 4\Omega \left(1 - \sqrt{\rho}\right)^{H_1} \sqrt{M} \Omega + 4\left(1 - \sqrt{\rho}\right)^{2H_1} M\Omega^2
	\\&\hspace{0.4cm} + 2 \Omega \left(1 - \sqrt{\rho}\right)^{H_1} \sqrt{M} \Omega + (2 + 24\gamma^2 L^2) \left(1 - \sqrt{\rho}\right)^{2H_1} M\Omega^2 + 24\gamma^2 L^2 \left(1 - \sqrt{\rho}\right)^{2H_1} \Omega^2   
	\\&\hspace{0.4cm} + \left(24\gamma^2 + \frac{32\gamma}{\mu}\right) L^2 \left(1 - \sqrt{\rho}\right)^{2H_1} \Omega^2 + \left(24\gamma^2 + \frac{32\gamma}{\mu}\right) L^2 \left(1 - \sqrt{\rho}\right)^{2H_1} M\Omega^2
	\\&\hspace{0.4cm} + 2 \Omega \left(1 - \sqrt{\rho}\right)^{H_1} \sqrt{M}\Omega
	+ 2 \left(1 - \sqrt{\rho}\right)^{2H_1} M\Omega^2
	\\&\hspace{0.4cm} + 2 \Omega \left(1 - \sqrt{\rho}\right)^{H_1} \sqrt{M}\Omega
	+ 2 \left(1 - \sqrt{\rho}\right)^{2H_1} M\Omega^2
	\\&\leq \left(48\gamma^2 + \frac{32\gamma}{\mu}\right) \cdot M(L\Omega + G)^2 \cdot (1 - \sqrt{\rho})^{2H_0}
	+ 10  \sqrt{M}\Omega^2 \cdot \left(1 - \sqrt{\rho}\right)^{H_1}
	\\&\hspace{0.4cm} 
	+\left(10 + 96\gamma^2 L^2 + \frac{64\gamma L^2}{\mu}\right) M\Omega^2 \cdot \left(1 - \sqrt{\rho}\right)^{2H_1}.
\end{align*}
Here we also use $\Omega$ and the same trick as \eqref{t101}. Then one can easy check that with our $H_0$ and $H_1$ from \eqref{H} it holds $\text{Err}(k) \leq \text{Err} \sim \varepsilon \mu \gamma$, then with \eqref{505} we get 
\begin{align*}
	\sqn{\bar z^{k+1} - z^*}
	&\leq 
	(1 - \gamma \mu)\sqn{\bar z^{k} - z^*} + \text{Err}.
\end{align*}
Running the recursion, we obtain
\begin{align*}
	\sqn{\bar z^{K} - z^*}
	&\leq 	(1 - \gamma \mu)^K\sqn{\bar z^{0} - z^*} + \frac{\text{Err}}{\gamma \mu},
\end{align*}
which completes the proof.
\EndProof

\textbf{Remark.} In the previous theorem, we obtained convergence along the point $\bar z^{K}$. This point is virtual and is not computed by the algorithm. But in fact, all local points $z^K_m$ are also very close to $\bar z^K$.

\textbf{Remark.} In this case \eqref{K3} dose not correspond to the number of communication rounds. To compute the number of rounds we need 
\begin{align*}
    K \times (H_0 + H_1) = \mathcal{\tilde O} \left( \frac{1}{\sqrt{\rho}}\left( 1 + \frac{\delta}{\mu} \right) \log \frac{\sqn{z^{0} - z^*}}{\varepsilon} \right).
\end{align*}
It is also easy to estimate the total number of local iterations on server:
\begin{align*}
    K \times T  &= 
    \cO \left( \frac{1}{\gamma \mu} \left(1 + \gamma L\right)\log\frac{1}{\tilde e} \log \frac{\sqn{z^{0} - z^*}}{\varepsilon} \right)
    \\&= \cO \left(  \left(\frac{1}{\gamma \mu} + \frac{L}{\mu}\right)\log\frac{1}{\tilde e} \log \frac{\sqn{z^{0} - z^*}}{\varepsilon} \right)
    \\&= \cO \left(  \left(1 + \frac{\delta}{\mu} + \frac{L}{\mu}\right)\log\frac{1}{\tilde e} \log \frac{\sqn{z^{0} - z^*}}{\varepsilon} \right).
\end{align*}


\subsubsection{Convex-Concave case}

This case is proved similarly to Theorem 6 (convergence) and Theorem 7 (inexact consensus). We just give the statement of the theorem:

\begin{theorem}
Let problem \eqref{inner} be solved by Extragradient with precision $e$:
\begin{align*}
    e = \mathcal{O}\left(\min \left\{ \frac{\varepsilon}{\delta} ; \frac{\varepsilon^2}{(L \Omega + G + \delta \Omega)^2} \right\} \right)
\end{align*}
and number of iterations $T$:
\begin{align*}
    T =  \cO  \left(\left(1 + \gamma L\right)\log\frac{\Omega^2}{e} \right).
\end{align*}
Suppose that parameters $H_0$ and $H_1$ satisfy
\begin{align*}
    H_0 = \mathcal{O} \left(\frac{1}{\sqrt{\rho}}\log \left(\frac{\left(\gamma^2 + \frac{\gamma}{\mu}\right) \cdot M(L\Omega + G)^2}{\varepsilon \gamma \mu}\right)\right), \notag\\ 
    H_1 = \mathcal{O} \left(\frac{1}{\sqrt{\rho}}\log \left(\frac{\left(1 + \gamma^2 L^2 + \frac{\gamma L^2}{\mu}\right) \cdot M\Omega^2}{\varepsilon \gamma \mu}\right)\right).
\end{align*}
Additionally, let us choose stepsize $\gamma$ as follows
\begin{align*}
    \gamma = \frac{1}{4\delta}.
\end{align*}
Then it holds that $\text{gap}(z^{K}_{avg}) \sim \varepsilon$ after 
\begin{align*}
    K = \cO \left( \frac{\delta \Omega^2}{ \varepsilon} \right)\quad \text{iterations},
\end{align*}
where $z^{K}_{avg}$ define as follows: $x^{K}_{avg} = \frac{1}{K}\sum_{k=0}^K u^{k}_x$, $y^{K}_{avg} = \frac{1}{K}\sum_{k=0}^K u^{k}_y$.
\end{theorem}

\section{Numerical Results}\label{app:numerical_results}

The numerical experiments are run on a machine with 8 Intel Core(TM) i7-9700KF 3.60GHz CPU cores with 64GB RAM. The methods are implemented in Python 3.7 using NumPy and SciPy.


In this section, we estimate the smoothness and strong convexity parameters for objectives used in all the experiments, as well as the similarity parameter. We denote the vector with all entries equal to one as $\one$ and the identity matrix as $I$ (with the sizes determined by the context). Given a set of data points $X = (x_1\ldots x_N)^\top\in\R^{N\times d}$ and an associated set of labels $y = (y_1\ldots y_N)^\top\in\R^N$, the Robust Linear Regression problem reads 
\begin{align*}
    \min_{\norm{w}\leq R_w} \max_{\norm{r}\leq R_r} g(w, r) := \frac{1}{2N} \sum\limits_{i=1}^N (w^T (x_i + r) - y_i)^2 + \frac{\lambda}{2} \| w\|^2 - \frac{\beta}{2} \|r \|^2.
\end{align*}
Note that we need constraints on $w$ to yield the bounds for smoothness and similarity parameters (this will be described below in this section). Equivalently, $g(w, r)$ can be expressed as
\begin{align*}
    g(w, r) = \frac{1}{2N} \norm{Xw + \one r^\top w - y}^2 + \frac{\lambda}{2}\sqn{w} - \frac{\beta}{2}\sqn{r},
\end{align*}
and its gradient w.r.t. $w$ and $r$ writes as
\begin{align*}
    \nabla_w g(w, r)
    &= \frac{1}{N} \cbraces{X^\top Xw + X^\top \one r^\top w - X^\top y + \one^\top (Xw - y)r} + rr^\top w + \lambda w, \\
    \nabla_r g(w, r)
    &= ww^\top r + \frac{1}{N} \one^\top (Xw - y)w - \beta r.
\end{align*}
The Hessian of $g(w, r)$ w.r.t. to $w$ and $r$ are
\begin{align*}
    \nabla_{ww}^2 g(w, r) &= \frac{1}{N} \cbraces{X^\top X + (X^\top\one r^\top + r\one^\top X)} + rr^\top + \lambda I, \\
    \nabla_{wr}^2 g(w, r) &= \frac{1}{N} \cbraces{X^\top\one w^\top + \one^\top (Xw - y) I} + r^\top w I + rw^\top, \\
    \nabla_{rr}^2 g(w, r) &= ww^\top - \beta I.
\end{align*}
We are now ready to estimate the spectrum of the Hessian taking into account the constraints on $w$ and $r$. For any $v\in\R^d$, we have
\begin{align*}
    \norm{\nabla_{ww}^2 g(w, r) v} 
    &\leq \frac{1}{N}\lambda_{\max}(X^\top X) \norm{v} + R_r^2 \norm{v} + \frac{1}{N} \norm{X^\top\one} R_r\norm{v} + \frac{1}{N} \norm{r\one^\top Xv} + \lambda\norm{v} \\
    &\leq \cbraces{\frac{1}{N} \lambda_{\max}(X^\top X) + R_r^2 + \frac{2}{N} R_r\norm{X^\top\one} + \lambda} \cdot \norm{v} =: L^g_{ww} \norm{v}, \\
    \norm{\nabla_{wr}^2 g(w, r) v}
    &\leq \frac{1}{N}\norm{X^\top\one w^\top r} + \frac{1}{N}\norm{\one^\top (Xw - y)v} + \norm{r^\top w v} + \norm{rw^\top v} \\
    &\leq \cbraces{\frac{2}{N}\norm{X^\top\one}R_w + \frac{1}{N}\one^\top y + 2 R_w R_r} \cdot \norm{v} =: L^g_{wr} \norm{v}, \\
    \norm{\nabla_{rr}^2 g(x, y) v} &\leq \norm{ww^\top v} + \beta\norm{v} \leq \cbraces{R_w^2 + \beta} \cdot \norm{v} =: L^g_{rr} \norm{v} .
\end{align*}
Therefore, we can estimate the Lipschitz constant of $\nabla g(w, r)$ as $L^g = \max(L^g_{ww}, L^g_{wr}, L^g_{rr})$.

Let us discuss the bound on the similarity parameter. Given two datasets $\braces{X\in\R^{N\times d},~ y\in\R^N}$ and $\braces{\widetilde X\in\R^{\widetilde N\times d},~ \widetilde y\in\R^{\widetilde N}}$, we define
\begin{align*}
    \widetilde g(w, r) = \frac{1}{2\widetilde N} \norm{\widetilde Xw + \one r^\top w - \tilde y}^2 + \frac{\lambda}{2}\sqn{w} - \frac{\beta}{2}\sqn{r}.
\end{align*}
To derive the similarity coefficient $\delta^{g, \widetilde g}$ between functions $g$ and $\widetilde g$, we separately estimate $\delta_{ww}^{g, \widetilde g},~ \delta_{wr}^{g, \widetilde g}$ and $\delta_{rr}^{g, \widetilde g}$.
\begin{align*}
    \delta^{g, \widetilde g}_{ww} &= \lambda_{\max}\cbraces{\frac{1}{N} X^\top X - \frac{1}{\widetilde N} \widetilde X^\top \widetilde X} + 2\norm{\frac{1}{N} X^\top\one - \frac{1}{\widetilde N} \widetilde
    X^\top\one} R_r, \\
    \delta^{g, \widetilde g}_{wr} &= 2\norm{\frac{1}{N} X^\top\one - \frac{1}{\widetilde N} \widetilde X^\top\one} R_w, \\
    \delta^{g, \widetilde g}_{rr} &= 0.
\end{align*}
We have $\delta^{g, \tilde g} = \max\{\delta^{g, \tilde g}_{ww}, \delta^{g, \tilde g}_{wr}\}$.

Finally, we estimate the strong convexity parameter as $\mu = \max(\lambda, \beta)$.

\end{document}